\documentclass[a4paper]{article}
\usepackage[T1]{fontenc}
\usepackage[english]{babel}
\usepackage[latin1]{inputenc}
\usepackage{dsfont}
\usepackage{graphicx}
\usepackage{color}
\usepackage{booktabs}
\usepackage{multirow}
\usepackage{subcaption} 
\usepackage{epsfig}
\usepackage{float}
\usepackage{cellspace}
\usepackage[normalem]{ulem}
\usepackage{diagbox}
\usepackage{amsthm,amssymb,amsbsy,amsmath,amsfonts,amssymb,amscd,mathrsfs}
\usepackage{colortbl}
\usepackage{verbatim}
\usepackage{url}
\usepackage{enumerate}
\usepackage{enumitem} %
\usepackage{hyperref} %
\usepackage{bbm}
\newtheorem{rem}{Remark}
\newtheorem{prop}{Proposition}[section]
\newtheorem{coro}{Corollary}[section]
\newtheorem{propri}{Property}[section]

\newtheorem{lem}{Lemma}[section]

\newtheorem{thm}{Theorem}[section]

\def\ds{\displaystyle}
\newcommand{\vphi}{\varphi}

\newcommand{\eps}{\varepsilon}
\newcommand{\R}{\mathbb{R}}

\newcommand\be{\begin{equation}}
\newcommand\ee{\end{equation}}
\numberwithin{equation}{section}
\newcommand{\N}{\mathbb{N}}
\providecommand{\norm}[1]{|#1|}

\title{Global stability of perturbed chemostat systems}
\author{C. Alvarez-Latuz\thanks{Avignon Universit\'e, Laboratoire de Math\'ematiques d'Avignon (EA 2151), F-84018 Avignon, France. {\tt claudia.alvarez-latuz@univ-avignon.fr}} , T. Bayen\thanks{Avignon Universit\'e, Laboratoire de Math\'ematiques d'Avignon (EA 2151), F-84018 Avignon, France. {\tt terence.bayen@univ-avignon.fr}} , J. Coville\thanks{UR 546 Biostatistique et Processus Spatiaux, INRAE, Domaine St Paul Site Agroparc, F-84000 Avignon, France.{\tt jerome.coville@inrae.fr}}}
\begin{document}
\maketitle

\begin{abstract}
This paper is devoted to the analysis of global stability of the chemostat system with a perturbation term 
representing any type of exchange between species. This conversion term depends on species and substrate concentrations but also on a positive perturbation parameter. 
After having written the invariant manifold as a union of a family of compact subsets, our main result states that  
for each subset in this family, there is a positive threshold for the perturbation parameter below which, the system is globally asymptotically stable in the corresponding subset. Our approach relies on the Malkin-Gorshin Theorem and on a Theorem by Smith and Waltman about perturbations of a globally stable steady state. 
Properties of steady-states and numerical simulations of the system's asymptotic behavior complete this study for two types of perturbation term between species. 

\end{abstract}

\section{Introduction}

{\bf{General context}}. The so-called {\it{chemostat system}} allows to model thanks to ordinary differential equations, the evolution of microbial species present in lakes or lagoons. In biotechnology, it can also be used to study the behavior of microbial species, with the aim, for example, of controlling the production of molecules of interest or reducing certain concentrations for water treatment \cite{Bastin90,dochain1}.
When the dilution rate is constant (this parameter allows to regulate the input substrate concentration) and in the presence of a single limiting substrate, the {\it{competitive exclusion principle}} (CEP) makes it possible to predict the asymptotical behavior of species concentrations. When considering non-decreasing kinetics for the species (such as Monod's kinetics \cite{Monod1,Monod2}), the CEP asserts that, generically, only one species survives asymptotically \cite{alain-livre,SmithWalt}. It is worth mentioning that in the chemostat system with a single limiting substrate, the species do not interact directly with each other but only indirectly through the substrate equation. 
\medskip

\noindent{\bf{Main objective}}. The aim of this paper is to study stability properties of the chemostat system when an exchange term\footnote{We will also use the terminology 
 ``perturbation term'' thoughout the paper.}  between  species is incorporated. This additional term in the system is quantified by a generally low perturbation parameter and  has several origins. It can represent, for example, mutation  between species (due to gene transfer) and in that case, the perturbation parameter can be seen as a mutation rate. To remain general, we prefer to use in this paper the terminology {\it{perturbed chemostat system}}, to describe the resulting dynamical system. 
In  presence of this phenomenon, each species is able to convert into neighboring species. 
Theoretical and numerical approaches to  study the chemostat system with a perturbation term as well as related questions can be found in \cite{Arkin,BCLC,BCM1,BM,BM2,deleen,deleen2,Coralie,lobry,szilard} (among others). In this paper, we are interested in addressing the next question : what can we say about the asymptotic behavior of a perturbed chemostat system particularly when the perturbation is ``small\footnote{This is made more precise in the rest of the paper.}''? Depending on the data defining the system  (such as dilution rate, perturbation parameter, kinetics), we wish to know if there is an invariant subset of the state space in which the perturbed system is globally asymptotically stable around a {\it{coexistence\footnote{This is an equilibrium point such that all species are present asymptotically.} steady-state}}.   
Even if the system in question is close to the chemostat system (in a way that has yet to be quantified), this question is not obvious. 
Indeed, it is well-known that the properties of a perturbed dynamical system can differ significantly from those of the non-perturbed system, potentially resulting in the emergence of multiple steady states or periodic orbits that would not exist without  perturbation. 

It is worth mentioning that for the chemostat system, previous studies (such as \cite{BCLC,deleen})  indicate that the techniques for studying the asymptotical  behavior of a perturbed version of the chemostat system differ from the classical approaches for the chemostat system (using for instance Lyapunov functions \cite{GMR2009,Hsu1,SmithWalt}).   
\medskip

\noindent{\bf{Related results}}. In order to highlight the novelty of our work, 
we would like first to recall some results in \cite{BCLC,deleen}. First, the article \cite{deleen} 
provides a first insight into the above question and it is as follows. Global stability is proven provided that  
kinetics are in a sufficiently narrow bundle (the limiting case being that they are all equal to a Monod type's kinetics) and that all yield coefficients are close to each other (the limiting case being that they are all equal to one). 
The paper \cite{BCLC} provides a similar global stability result provided that kinetics are of Monod type and that the dilution rate is small enough
 (the limiting case being that the dilution rate is zero). Additionally, the occurrence of a single coexistence steady-state (apart from the washout) is established which differs from the chemostat system where several steady states occur (the number of them being equal to the number of species apart from the wash-out). 
This work is in line with the preceding studies, but, it is probably set in the most natural case, {\it{i.e.}}, when the kinetics are monotone (but arbitrarily),
 when the dilution rate is also arbitrarily (which is more realistic than in \cite{BCLC} from an experimental point of view) and  
when the perturbation parameter tends to zero. 
Furthermore, we also consider  
the case of a non-linear exchange term depending both on the substrate and species concentrations 
and non-necessarily linear w.r.t.~species concentrations as in \cite{BCLC}. 
\medskip

\noindent {\bf{Main contribution}}.  
In order to study global stability of a general perturbed dynamical system (knowing that the non perturbed system is globally asymptotically stable), Smith and Waltman introduced in \cite[Theorem 2.2]{SW99} (see also their Corollary 2.3) sufficient conditions on the system to guarantee that global stability is preserved as the perturbation parameter tends to zero. In order to apply this fundamental result in our setting, it is essential to verify one of the hypotheses which states that all trajectories of the perturbed system should enter after a certain time into a (fixed) compact subset of the invariant manifold and remain in it. In our setting, this hypothesis is not evident to verify. A good approach is to use the Malkin-Gorshin Theorem \cite{tewfik0} (see also \cite{tewfik1,tewfik2}). This result is part of the folklore and it allows to show that trajectories of the perturbed system remain close to trajectories of the non-perturbed one provided that the perturbation parameter is small enough. It is proved in \cite{tewfik0} in the context of B-stability which in some sense is a weaker notion than local stability. The B-stability of a differential system is verified if the original system is exponentially stable. 
Since the chemostat system fulfills the latter condition, we shall not enter into the notion of B-stability throughout this paper. In order to combine the Malkin-Gorshin Theorem and \cite[Theorem 2.2]{SW99}, we need to introduce a family of compact subsets of the invariant manifold (to apply the  Malkin-Gorshin Theorem uniformly w.r.t.~the initial condition). Our main result states that for each subset in this family, there is a positive threshold for the perturbation parameter below which, the system is globally asymptotically stable in the corresponding subset (Theorem~\ref{main1}). We complete this study by analyzing steady-states when the perturbation term depends on the substrate but is linear w.r.t.~concentration species. Numerical simulations are carried out for two types of exchange terms and highlight the behavior of the perturbed dynamical system as expected thanks to Theorem~\ref{main1}. 
\medskip

\noindent{\bf{Organization of the paper}}. The paper is structured as follows. In Section \ref{sec:pre}, we present the model and our main hypotheses. Additionally, we provide Theorem~\ref{prop-claudia-ext-SW} that is a slight extension of \cite[Corollary 2.3]{SW99} to our setting and Theorem~\ref{prop-claudia-MG} which provides a simplified version of the Malkin-Gorshin Theorem. For sake of completeness, these propositions are proved in the Appendix. In section~\ref{sec:3}, we analyze in details steady-states of the perturbed system in the case where the perturbation term is linear w.r.t.~species concentrations. Finally, Section~\ref{sec:simu} develops numerical simulations of the perturbed chemostat system for two types of exchange terms. Section~\ref{sec:conc} explores various perspectives derived from the results obtained in this paper.

\section{Model presentation and main result}
\label{sec:pre}
\medskip

\subsection{The chemostat system with a general perturbation term}

Throughout this paper, we consider the chemostat system with a single limiting substrate and $n \geq 1$ species, and we also incorporate an exchange term between them that depends on species concentrations, substrate concentration, and on a perturbation parameter called $\eps$. This yields the following dynamical system
\begin{equation}{\label{sys1}}
\left\{
\begin{array}{cl}
\dot{x}&=D(s)x-u x + h(x,s,\eps), \smallskip\\
\dot{s}&=\ds -\sum_{j=1}^n \frac{\mu_j(s)x_j}{Y_j}+u(s_{in}-s),
\end{array}
\right.
\end{equation}
where: 
\begin{itemize}
\item[$\bullet$] $x:=(x_1,...,x_n)^\top\in\R^n$ is a vector containing the species concentrations $x_i$ for  $i\in \{1,...,n\}$ ;  \\
\item[$\bullet$] $s$ and $s_{in}$ denote respectively the substrate and input substrate concentration ($s_{in}>0$ is fixed) ;
\item[$\bullet$] $D(s)\in \R^{n\times n}$ denotes the diagonal matrix $\mathrm{diag}(\mu_1(s),...,\mu_n(s))$ where for all $i\in \{1,...,n\}$, the function $\mu_i:\R_+\to \R_+$ is the kinetics of species $i$ ; 
\item[$\bullet$] for all $i\in \{1,...,n\}$, $Y_i>0$ is the yield coefficient associated with species $i$ ;
\item[$\bullet$] $u\geq 0$ is the dilution rate and $\eps\geq 0$ is the perturbation parameter ; 
\item[$\bullet$] the function $h:\R^n\times\R\times\R\to \R^n$ represents the perturbation term so that for all $i\in \{1,...,n\}$, $h_i(x,s,\eps)$ describes how species $i$ converts into other species at a rate $\eps$. 
\end{itemize}

The perturbation parameter $\eps$ plays different roles depending on the interaction type considered between species (typically, it can represent a mutation rate or a constant probability of mutation). Regardless of the biological meaning of $\eps$, it is relevant in our formulation as it will allow us to see \eqref{sys1} as a perturbation of the chemostat system. Let us now introduce the following hypotheses on the data. Hereafter, $|\xi|$ stands for the euclidean norm of a vector $\xi\in \R^n$. 

\begin{enumerate}[label={\rm(A\arabic*)}]
    \item\label{as:1} For all $i\in \{1,...,n\}$,  $\mu_i$ is of class $C^1$, non-decreasing, bounded over $\R_+$, and such that $\mu_i(0)=0$. 
  \item\label{as:2} The function $h$ is of class $C^1$ and satisfies  $h(x,s,0)= 0$ for all $(x,s)\in \R^{n+1}$. Moreover, it is  with linear growth, {\it{i.e.}},  for all $\eps>0$ there is $c_\eps \geq 0$ such that for all $(x,s)\in \R^{n+1}$, 
  \begin{equation}{\label{LG}}
    |h(x,s,\eps)| \leq c_\eps\big(|x|+1\big).
    \end{equation}
       \item\label{as:3} For all $i\in \{1,...,n\}$ and all $x\in \R^n$, if $x_i=0$, then for all $(\eps,s)\in \R_+^*\times [0,s_{in}]$, one has $h_i(x,s,\eps)\geq 0$. 
    \item\label{as:4} For every $(x,s,\eps)\in \R^n\times\R\times \R_+$, one has
     
    $$\displaystyle\sum_{j=1}^n h_j(x,s,\eps)=0.$$
    
\end{enumerate}

Some comments on the preceding hypotheses are in order. 
\begin{rem}{\label{hypo-remark}}
{\rm{(i)}}. Hypothesis \ref{as:1} is standard when considering  chemostat type's systems (see, {\it{e.g.}}, {\rm{\cite{SmithWalt}}}). However, depending on the application model, other types of kinetics are relevant such as Haldane\footnote{Such kinetics allow to take into account inhibition through substrate.} type's kinetics which are non-monotonic or kinetics with a more complicated behavior such as in {\rm{\cite{BGM13}}}. In this paper, we shall restrict our attention to monotone kinetics only.  
\smallskip
\\
{\rm{(ii)}}. Supposing that $h$ vanishes for $\eps=0$ (see \ref{as:2}) allows us to retrieve the chemostat system 
\begin{equation}{\label{sys2}}
\left\{
\begin{array}{cl}
\dot{x}_i&=\mu_i(s)x_i-u x_i, \quad 1 \leq i \leq n \smallskip\\
\dot{s}&=\ds -\sum_{j=1}^n \frac{\mu_j(s)x_j}{Y_j}+u(s_{in}-s),
\end{array}
\right.
\end{equation}
for $\eps=0$. This hypothesis is important in our study in order to investigate global stability properties of \eqref{sys1}. Note that \eqref{LG} is a standard hypothesis in the theory of ordinary differential equations in order to prevent the blow-up phenomenon of solutions. 
\smallskip
\\
{\rm{(iii)}}. Hypothesis \ref{as:3} is  essential for biological purpose to guarantee that all species remain non-negative over time  
and Hypothesis \ref{as:4} implies that mass is conserved during the exchange process between species (this property is crucial in our study, particularly for proving Lemma~\ref{lem0} and Proposition~\ref{lem:invariance}).  
\end{rem}

\subsection{Invariant domain and related properties}
In this section, we introduce the invariant domain 
that will be considered throughout the paper.  
The next lemma is a preliminary result in order to introduce this set. 

\begin{lem}{\label{lem0}}
    Suppose that hypotheses \ref{as:1}-\ref{as:2}-\ref{as:3} are satisfied. 
\smallskip
\\
{\rm{(i)}} Then, for every $(x_0,s_0)\in \R^n_+\times[0,s_{in}]$ and for every $(\eps,u)\in \R_+\times \R_+^*$, there is a unique solution of \eqref{sys1} defined over $\R_+$. Moreover,  $\R^n_+\times[0,s_{in}]$ is forward invariant. 
\smallskip\\
{\rm{(ii)}} For every $(s,\eps)\in [0,s_{in}]\times \R_+$, one has $h(0,s,\eps)=0$ and the point  $E_{wo}:=(0_{\R^n},s_{in})$ is an equilibrium of \eqref{sys1} called the wash-out steady-state. 
\end{lem}
\begin{proof} 
The existence of a unique solution to \eqref{sys1} essentially follows from \eqref{LG} which guarantees that solution are defined globally over $\R_+$. The fact that 
$\R^n_+\times[0,s_{in}]$ is invariant follows from \ref{as:3} and using that $\dot{s}_{|_{s=0}}=us_{in} \geq 0$ and 
$\dot{s}_{|_{s=s_{in}}}\leq 0$. This takes care of (i). To prove (ii), observe that if $x=0$, then for every $1 \leq i \leq n$, one has 
$h_i(0,s,\eps)\geq 0$ thanks to \ref{as:3}. Now, \ref{as:4} implies that $\sum_{j=1}^n h_j(0,s,\eps)=0$, hence, we must have $h_i(0,s,\eps)=0$ for every $1 \leq i \leq n$ as wanted. The fact that $E_{wo}$ is a steady-state of \eqref{sys1} follows from the equality $h(0,s_{in},\eps)=0$ for every $\eps\geq 0$. This ends the proof. 
\end{proof}
From Lemma~\ref{lem0}, $[0,s_{in}]$ is forward invariant for $\dot{s}$ (and also attractive in $\R_+\backslash [0,s_{in}]$), hence, without any loss of generality, we will suppose that $s(0)\in [0,s_{in}]$, throughout the paper.  
In addition, it is readily seen that $s(t)>0$ for all time $t>0$ (even if $s(0)=0$) provided that $u>0$. 
Regarding the positivity of the variables $x_i$, obtaining such a property depends on the assumptions made about the function $h$.
The next lemma provides a sufficient condition ensuring positivity of $x_i(t)$ for $t>0$ and $1 \leq i \leq n$.     
\begin{lem}{\label{lem-positivity}}
Suppose that \ref{as:1}-\ref{as:2}-\ref{as:3} are satisfied and that $h$ satisfies the  additional condition : for every triplet
$(x,s,\eps)\in \R_+^n\times[0,s_{in}] \times \R_+$ and for every index $i\in \{1,...,n\}$: 
\begin{equation}{\label{CS-pos}}
\left\{\begin{array}{clcccl}
\mathrm{if} & i \notin \{1,n\} & \mathrm{and}   & x_i=h_i(x,s,\eps)=0 & \Rightarrow & x_{i-1}=x_{i+1}=0,\\
\mathrm{if} & i=1 & \mathrm{and}  & x_i=h_i(x,s,\eps)=0 & \Rightarrow & x_{i+1}=0,\\
\mathrm{if} & i=n & \mathrm{and} & x_i=h_i(x,s,\eps)=0 & \Rightarrow & x_{i-1}=0.
\end{array}\right.
\end{equation}
Then, for every initial condition $(x_0,s_0)\in \R_+^n \times [0,s_{in}]$ such that $x_0\not=0$, the unique solution to \eqref{sys1} with initial condition $(x_0,s_0)$ at time $t=0$ satisfies $x_i(t)>0$ for every $1 \leq i \leq n$.
\end{lem}
\begin{proof}
Take $1 \leq i\leq n$ such that $x_i(0)>0$ and suppose by contradiction that $x_i(\cdot)$ vanishes over $\R_+$. Let then $t_0:=\inf\{t>0 \; ; \; x_i(t)=0\}$. 
Since $x_i>0$ over $[0,t_0)$, we deduce that $\dot{x}_i(t_0)\leq 0$ so that 
$h_i(x(t_0),s(t_0),\eps)=x_i(t_0)=0$. If $1<i<n$, we  get
$x_{i-1}(t_0)=x_{i+1}(t_0)=0$. By a similar argumentation we deduce that 
$h_{i-1}(x(t_0),s(t_0),\eps)=h_{i+1}(x(t_0),s(t_0),\eps)=0$ so that if $i-1 \geq 1$ and $i+1 \leq n$, condition \eqref{CS-pos} implies that $x_{i-2}(t_0)=x_{i+2}(t_0)=0$. 
By an inductive reasoning, we deduce that $x(t_0)=0$, hence, Lemma~\ref{lem0} implies that $x(t)=0$ for every $t\geq 0$ which is a contradiction since $x(0)\not=0$. 
Hence, we deduce that $x_i(t)>0$ for every $t>0$. The same reasoning applies to all other indexes as well which ends the proof.  
\end{proof}
We give in Section~\ref{sec:simu} two examples where the latter hypothesis about $h$ is fulfilled implying that every species is present at any time $t>0$ (in contrast with the chemostat system where for all $1 \leq i \leq n$, either $x_i\equiv 0$ if $x_i(0)=0$ or $x_i(t)>0$ for all time $t\geq 0$). 


We now turn to the definition of the invariant set $\Delta$ for \eqref{sys1}. Set $Y_+ := \max~\big\{ Y_1,\ldots,Y_n\big\}$ and let $\Delta$ be defined as
\begin{equation}
\Delta:=\bigg{\{}(x,s)\in \R_+^n \times [0,s_{in}]\; ; \; s+ \frac{1}{Y_+}\sum_{j=1}^n x_j \leq s_{in}\bigg{\}}. 
\end{equation}
This set will play a significant role in the remainder of the paper (in particular to state our global stability results). 
Hereafter, we say that a set is (forward) invariant for a dynamical system if every solution starting in this set at time $t=0$ 
remains in it for all time $t\geq 0$. 
We say that it is attracting if the distance between this set and solutions starting outside the set at time $t=0$ 
converges to zero as $t\rightarrow +\infty$.

\begin{prop}\label{lem:invariance}
If hypotheses \ref{as:1}-\ref{as:2}-\ref{as:3}-\ref{as:4} are satisfied, then $\Delta$ is an invariant and attractive manifold for solutions of \eqref{sys1}. Moreover, for every  $(x_0,s_0)\in \R_+^n \times [0,s_{in}]$, 
either the unique solution to \eqref{sys1} such that $(x(0),s(0))=(x_0,s_0)$ satisfies $b(t)\rightarrow s_{in}$ as $t\rightarrow +\infty$ or it enters into $\Delta$ in finite time.  
\end{prop}
\begin{proof}
Let us set $b:=s+\frac{1}{Y_+}\sum_{j=1}^n x_j$. We have 
$$
\dot{b}  = \sum_{j=1}^n\left(\frac{1}{Y_+}-\frac{1}{Y_j}\right) \mu_j(s)x_j+u\left( s_{in}-b \right)\leq u(s_{in}-b),
$$
from which we deduce that for every time $t\geq 0$, one has
\begin{equation}{\label{tmp-propri}}
b(t) \leq s_{in}+\big( b(0)-s_{in}\big)e^{-u t},
\end{equation}
implying that $\Delta$ is invariant through \eqref{sys1}. Take an initial condition $(x(0),s(0))\in (\R_+^n \times [0,s_{in}])\backslash \Delta$ and 
let $(x(\cdot),s(\cdot))$ be the unique solution to \eqref{sys1} starting at $t=0$ from this initial condition. If there is a time $t'> 0$ such that 
$b(t') \leq s_{in}$, we obtain using a similar argumentation as \eqref{tmp-propri} that $b(t)\leq s_{in}$ for all time $t\geq t'$ which shows that the distance of the solution to $\Delta$ is zero for all $t\geq t'$. On the other hand, if 
for every $t\geq 0$, one has $b(t)>s_{in}$, we obtain that 
$$
s_{in} < b(t)  \leq s_{in}+\big( b(0)-s_{in}\big)e^{-u t},
$$
for every $t \geq 0$. 
Hence, the distance of the solution to $\Delta$ (that is proportional to 
$ b(t) -s_{in}$ since the active part of boundary of $\Delta$ is the hyperplane of equation $s+\frac{1}{Y_+}\sum_{j=1}^n x_j=s_{in}$ in $\R_+^n \times [0,s_{in}]$) goes to zero as $t\rightarrow +\infty$. To prove the rest of the proposition, suppose that $t\mapsto b(t)$ does not converge to $s_{in}$ as $t\rightarrow +\infty$. 
If $b(t)>s_{in}$ for all time $t\geq 0$, the function $b(\cdot)$ necessarily converges to $s_{in}$ from \eqref{tmp-propri}. 
Hence, there is $t_0\geq 0$ such that $b(t_0)\leq s_{in}$. 
This proves exactly that the trajectory has entered into $\Delta$ in finite time. 
\end{proof}
\begin{rem} 
{\rm{(i)}}. In Lemma~\ref{lem0} and in Proposition~\ref{lem:invariance}, \ref{as:1} can be weakened supposing that $\mu_i$ is only non-negative, continuous and null at $s=0$, but, as said in Remark~\ref{hypo-remark}, we shall not consider non-monotonic kinetics in this paper. 
\smallskip
\\
{\rm{(ii)}} In the case of the chemostat system, it is relevant to introduce $b_0:=s+\sum_{j=1}^n x_j$ so that changing $x_i$ into $x_i/Y_i$ if necessary, $b_0$ satisfies 
$\dot{b}_0=u(s_{in}-b_0)$. In that case, $\Delta_0:=\{(x,s) \in \R_+^n \times [0,s_{in}]\; ; \; \sum_{j=1}^n x_j+s=s_{in}\}$ is an invariant and attractive manifold of \eqref{sys2}, but, solutions to \eqref{sys2} starting outside this set never enter into $\Delta_0$ in finite time. 
A key difference that arises whenever $\eps>0$ is that solutions to \eqref{sys1} with $\eps>0$ may enter the region $\Delta$ in finite time, provided that the function $b$ does not converge to $s_{in}$. 
We highlight this property numerically in Section~\ref{sec:simu} (see Fig.~\ref{fig:monod_lin2} and Fig.~\ref{fig:monod_cul2}). 
\end{rem}

\subsection{Competitive exclusion principle}

\begin{figure}[H]
    \centering
    \includegraphics[width=0.4\linewidth]{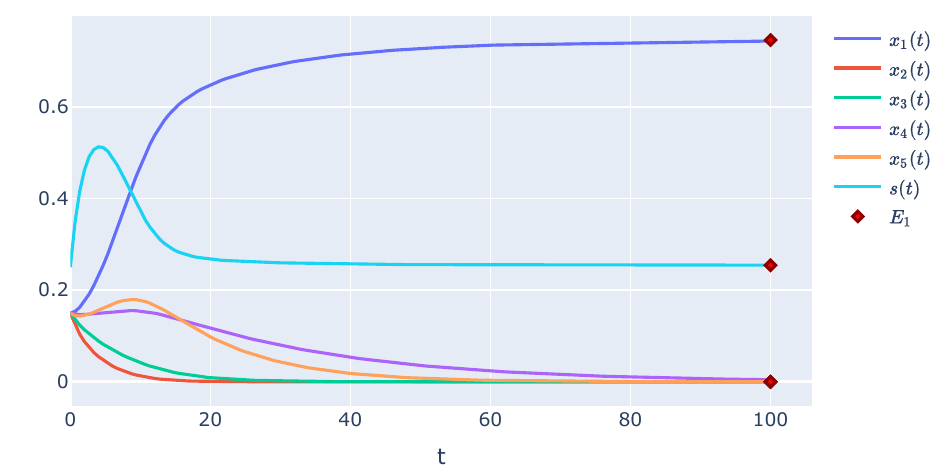}
    \includegraphics[width=0.4\linewidth]{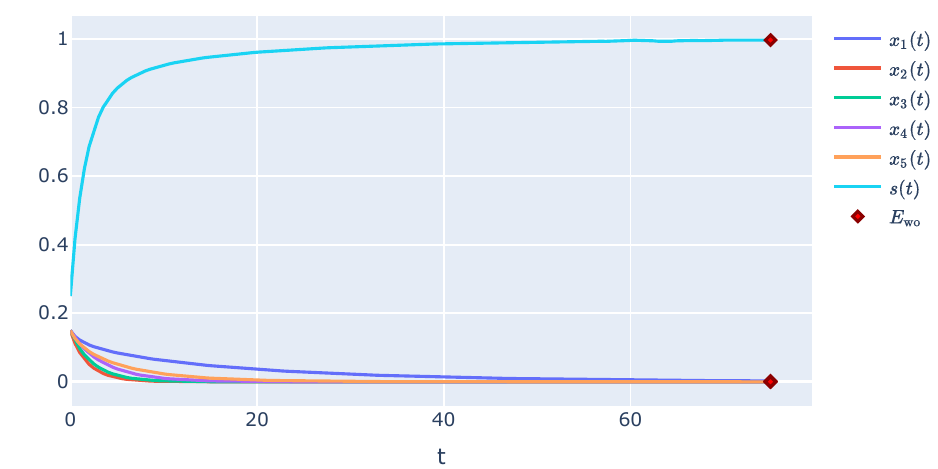}
    \caption{Illustration of Theorem~\ref{CEP-thm} for {\small{$(x_0,s_0)=(0.15,0.15,0.15,0.15,0.15,0.25)$}}. {\it{Fig.}}~left: 
    convergence to {\small{$E_1=(0.75,0,0,0,0,0.25)$}} with $u=0.4$ ; {\it{Fig.}}~right: convergence to {\small{$E_{wo}=(0,0,0,0,0,1)$}} with $u=0.7$.}
        \label{fig:CEP}
\end{figure}

The analysis of stability properties of \eqref{sys1} when the perturbation parameter tends to zero is closely related to properties of the chemostat system \eqref{sys2} corresponding to \eqref{sys1} for $\eps=0$, that is why, we wish now to recall the competitive exclusion principle (see, {\it{e.g.}}, \cite{alain-livre,SmithWalt}). In what follows, we set
$$
\vphi(u):=\min \big\{ \mu^{-1}_i(u)\; ; \; i=1,...,n\big\} \in [0,+\infty]. 
$$
Note that $\vphi(u)$ can be equal to $+\infty$ if $u$ is greater than $\max_{s\in [0,s_{in}]} \mu_i(s)$ for every $1 \leq i \leq n$. 
For $\eps=0$, \eqref{sys1} has at most exactly $n+1$ steady-states (depending on the value of $u$), namely, 
$$
E_{i}:=(0,...,0,1-\mu^{-1}_i(u),0,...,0,\mu^{-1}_i(u))\in \R^{n+1}, \quad 1 \leq i \leq n,
$$
together with $E_{wo}$. 
The CEP can be formulated as follows. 
\begin{thm}
{\label{CEP-thm}} Suppose that \ref{as:1} is fulfilled.  
If there is a unique $i^*\in\{1,\ldots,n\}$ such that $\mu^{-1}_{i^*}(u)=\vphi(u)<+\infty$, 
then, for every initial condition $(x_0,s_0)\in \R_+^n \times [0,s_{in}]$ such that $x_{i^*}(0)>0$, the unique solution of \eqref{sys1} for $\eps=0$ starting at $(x_0,s_0)$ at time $0$, converges to $(x^*,s^*)=E_{i^*}$, and this steady state is locally exponentially stable. If $\vphi(u) =+\infty$, then,  for every initial condition $(x_0,s_0)\in \R_+^n \times [0,1]$, the unique solution to \eqref{sys1} with $\eps=0$ starting at $(x_0,s_0)$ at time $0$ converges to 
$E_{wo}$. \end{thm}
Fig.~\ref{fig:CEP} illustrates Theorem~\ref{CEP-thm} in two cases (kinetics are given in Table~\ref{tab:parameters}) with five species and $s_{in}=1$. The first case is with $u=0.4$ where species $1$ dominates the culture ($\vphi(u)=0.25$) whereas case 2 is with $u=0.7$ so that wash-out occurs leading to extinction of species ($\vphi(u)=+\infty$).  
As it is well-known the behavior of a perturbed dynamical system may slightly differ from the behavior of the non-perturbed system. In our setting and as far as we know, the available methods to prove the CEP (such as  Lyapunov functions) cannot be straightforwardly transferred to study \eqref{sys1}. 
That is why, we will enunciate in the next section intermediate results before stating the main theorem about the behavior of \eqref{sys1}. 

\subsection{Global stability of the perturbed chemostat system}
In this section, we prove Theorem \ref{main1} which is our main result. We start by recalling a result by Smith and Waltman from \cite{SW99}, but adapted to our context (Theorem~\ref{prop-claudia-ext-SW}) and we also state a weaker version of the Malkin-Gorshin Theorem\footnote{This result is demonstrated in \cite{tewfik0} in the context of B-stability. However, we deal here with asymptotic stability which is a stronger notion but sufficient for our purpose.} to be found in \cite{tewfik0} (Theorem~\ref{prop-claudia-MG}).  Doing so, let us consider a general Cauchy problem
\begin{equation}{\label{eq:generalode}}
\left\{
\begin{array}{cl}
\dot{y}&= f(y,\eps), \\
y(0)&=y_0,
\end{array}
\right.\end{equation}
where $f:\mathcal{V} \times \R \rightarrow \R^m$ is of class $C^1$, $\mathcal{V}\subset \R^m$ is such that $\mathrm{Int}(\mathcal{V})\not=\emptyset$\footnote{Throughout the paper, the interior of a set $A\subset \R^m$ is denoted by $\mathrm{Int}(A)$.}, $y_0\in \mathcal{V}$, and $\eps\in \R_+$ is a parameter. 
We suppose that for every $y_0\in \R^m$, the unique solution to \eqref{eq:generalode} denoted by $y(\cdot,y_0,\eps)$ is defined over $\R_+$ and that $\mathcal{V}$ is forward invariant by \eqref{eq:generalode}. The next result\footnote{We raise this result
 as a theorem since \cite[Corollary 2.3]{SW99} is just a consequence of \cite[Theorem 2.2]{SW99} which is one of the main result of \cite{SW99}.} is a slight adaptation of \cite[Corollary 2.3]{SW99}. 
\begin{thm}{\label{prop-claudia-ext-SW}} 
Let $\mathcal{U}$ be a subset of $\R^m$ such that $\mathcal{U}\subset\mathcal{V}$ and $\mathrm{Int}(\mathcal{U})\neq\emptyset$.  Suppose that there is 
$y^*\in \mathrm{Int}(\mathcal{U})$ such that $f(y^*,0)=0$, $D_yf(y^*,0)$ is Hurwitz, and $y^*$ is attracting\footnote{This means that $y(t,y_0,0)\to y^*$ for every $y_0\in \mathcal{U}$ as $t\rightarrow+\infty$.} for solutions of \eqref{eq:generalode} in $\mathcal{U}$ with $\eps=0$. Moreover, assume that the following hypothesis is fulfilled:
    \begin{enumerate}[label={\rm(H\arabic*)}]
        \item\label{hyp:1} There exist a compact set $\mathcal{C}\subset \mathcal{U}$ and $\bar{\eps}>0$ such that for each $\eps\in[0,\bar{\eps}]$ and for each $y_0\in \mathcal{U}$, the unique solution to \eqref{eq:generalode} 
        reaches $\mathcal{C}$ meaning that $y(t,y_0,\eps)\in \mathcal{C}$ for every $t$ large enough.
    \end{enumerate}
Then, there are $\bar{\eps}_0\in(0,\bar{\eps})$ and a unique point $\hat{y}(\eps)\in \mathcal{U}$ such that for every $\eps\in [0,\bar{\eps}_0]$, one has 
$f(\hat{y}(\eps),\eps)=0$ and  $y(t,y_0,\eps) \rightarrow \hat{y}(\eps)$ for all $y_0\in \mathcal{U}$ as $t\rightarrow +\infty$. 
\end{thm}
The only difference with \cite[Corollary 2.3]{SW99} is that $\mathcal{U}$ is not supposed to be invariant but rather contained in another 
set $\mathcal{V}$ which is invariant. For this reason, we propose a proof of this proposition in the Appendix following~\cite{SW99} 
(also, the proof is made in the context of ordinary differential equations and not in an abstract setting).    
Finally, we would like to insist on the fact that the point $y^*$ may belong to the boundary of the set 
$\mathcal{U}$ if the dynamics $f$ can be extended to an open set around $y^*$ 
whose intersection with $\mathcal{U}$ is convex (as it is mentioned in \cite[Remark 2.1]{SW99}).

We now turn to the Malkin-Gorshin Theorem \cite[Theorem 4]{tewfik0} which is part of the folklore in the theory of dynamical systems involving regular perturbations. A  version of this result involving the concept of B-stability can be found in \cite{tewfik0}. 
We provide below a version of this result in the  (stronger) setting of asymptotic stability which is enough for our purpose.

\begin{thm}{\label{prop-claudia-MG}}
Let $y^*\in\R^m$ be such that $f(y^*,0)=0$. 
Suppose that $D_yf(y^*,0)$ is Hurwitz and that $y^*$ is attracting for solutions of \eqref{eq:generalode} with $\eps=0$ and initial conditions in a non-empty compact set $\mathcal{C}\subset\R^m$. Then, for every $\delta>0$, there is $\bar{\eps}>0$ such that for all $\eps\in[0,\bar{\eps}]$ and all $y_0 \in \mathcal{C}$, 
    \[
    \sup_{t\geq 0}\ |y(t,y_0,\eps)-y(t,y_0,0)| \leq \delta. 
    \]
\end{thm}
The proof of this result (using local stability) can be found in the Appendix. 
Let us emphasize the fact that the parameter $\delta$ depends on a compact set $\mathcal{C}$ for initial conditions (in contrast with \cite[Theorem 4]{tewfik0} which involves a single initial condition). In our approach, considering such a compact set will be needed to prove Theorem~\ref{main1} 
whose proof relies on the combination of Theorem~\ref{prop-claudia-ext-SW} and Theorem~\ref{prop-claudia-MG} in the setting of \eqref{sys1}. More precisely, the result of Theorem~\ref{prop-claudia-MG} will allow us to verify \ref{hyp:1} in Theorem~\ref{prop-claudia-ext-SW}.

Coming back to \eqref{sys1}, 
for every $1 \leq i \leq n$ and for every $\alpha\in (0,s_{in})$, we define a non-empty compact subset $\Delta_{i,\alpha}$ of $\R^n$ as

$$
\Delta_{i,\alpha}:=\{(x,s)\in \Delta \; ; \; x_{i}\geq \alpha\},
$$ 

so that for a fixed index $1 \leq i \leq n$, the union of these sets is associated with $\Delta$ through the following relationship:
$$
\Delta\backslash \{(x,s)\in \R_+^n \times [0,s_{in}] \; ; \; x_i=0\} = \bigcup_{\alpha>0} \Delta_{i,\alpha}. 
$$
Note that the sets $\Delta_{i,\alpha}$ are not necessarily invariant for \eqref{sys1}. 
Our main result states that \eqref{sys1} has a steady-state (apart from the wash-out) that 
 is  attracting in subsets $\Delta_{i^*,\alpha}$ provided that the perturbation parameter is small enough 
 (where given a value of the dilution rate $u$, $i^*$ represents the index of the species that dominates the competition in Theorem~\ref{CEP-thm}).  

\begin{thm}{\label{main1}}
Suppose that hypotheses \ref{as:1}-\ref{as:2}-\ref{as:3}-\ref{as:4} are fulfilled. Let $u>0$ be such that there is a unique $i^*\in\{1,...,n\}$ such that $\mu^{-1}_{i^*}(u)=\vphi(u)<+\infty$ and let $\alpha\in(0,s_{in}-\vphi(u))$. Then, there is $\eps_{u,\alpha}>0$ such that for every $\eps\in [0,\eps_{u,\alpha}]$, 
\eqref{sys1} has a steady-state $(x^{\eps,u},s^{\eps,u})\in\Delta_{i^*,\alpha}$ apart from the washout.
Moreover, for every $(x_0,s_0)\in \Delta_{i^*,\alpha}$, the unique solution to \eqref{sys1} starting at $(x_0,s_0)$ at time $t=0$, converges to  $(x^{\eps,u},s^{\eps,u})$. 
\end{thm}
\begin{proof}
Let  $f:\R^{n}_+\times [0,s_{in}] \times \R_+\rightarrow\R^{n+1}$ be defined as  
\begin{equation}{\label{def-f}}
f(y,\eps):= \left( \begin{array}{cccc}
 (\mu_1(s)-u)x_1 + h_1(x,s,\eps)    \\
 \vdots    \\
  (\mu_n(s)-u)x_n + h_n(x,s,\eps)     \vspace{0.5em} \\ 
    -\sum_{j=1}^n  \frac{\mu_j(s)x_j}{Y_j}+u(s_{in}-s)
\end{array}\right)
\end{equation}
for every $y:=(x,s)\in \R_+^{n}\times [0,s_{in}]$ and all $\eps \geq 0$. Set $y^*:=(x^*,s^*)$ where $(x^*,s^*)=E_{i^*}$ so that $f(y^*,0)=0$ and recall that $y_{i^*}^*=s_{in}-\vphi(u)$, hence, 
$0<\alpha<y_{i^*}^*$. For every $y_0:=(x_0,s_0)\in \R_+^n \times [0,s_{in}]$, we denote by $y(\cdot,y_0,\eps)$ the unique solution to \eqref{sys1} starting at $y_0$ at time $t=0$.
In virtue of Theorem~\ref{CEP-thm}, we have $y(t,y_0,0)\rightarrow y^*$ for every $y_0\in\R_+^{n}\times[0,s_{in}]$ whose $i^*$-th component is positive at time $t=0$. It follows that there is an instant $t_0\geq 0$ (that depends on the initial condition $y_0$) such that 
$$y_{i^*}(t,y_0,0)\geq y_{i^*}^*-\eta,
$$ for every $t\geq t_0$ where $\eta>0$ is such that $y_{i^*}^*-\eta>\alpha$ ($\eta$ does not depend on $y_0$).  
Now, the Jacobian matrix of $f$ w.r.t.~$y$ at $(y^*,0)$ is just the Jacobian matrix of the chemostat system at the globally asymptotically stable steady-state $E_{i^*}$. 
Hence, it is of Hurwitz type (this property is standard, see, {\it{e.g.}}, \cite{alain-livre,SmithWalt}). 
Moreover, for $\eps=0$, the CEP implies that $y^*$ is attracting for solutions to \eqref{sys1} (with $\eps=0$) 
and starting in $\Delta_{i^*,\alpha}$ at time $t=0$. 
We are now in a position to apply Theorem~\ref{prop-claudia-MG} to $f$ at the point $(y^*,0)$ with $\mathcal{C}=\Delta_{i^*,\alpha}$.  
Take $\delta>0$ such that $\delta<y_{i^*}^*-\eta-\alpha$. 
It follows that there is $\bar{\eps}>0$ such that for every $\eps\in[0,\bar{\eps}]$, for every $y_0 \in \mathcal{C}$, and for every $t\geq 0$, 
$$|y_{i^*}(t,y_0,0)-y_{i^*}(t,y_0,\eps)|\leq \delta.$$
Now, take $y_0\in \mathcal{C}$. We deduce that there is $t_0 \geq 0$ (as before, depending on $y_0$) such that for every $t\geq t_0$, 
\begin{align*}
    y_{i^*}(t,y_0,\eps)\geq y_{i^*}(t,y_0,0) - \delta \geq y_{i^*}^*-\eta-\delta\geq\alpha. 
\end{align*}
Therefore, $y(t,y_0,\eps)\in \Delta_{i^*,\alpha}$ for every $\eps\in[0,\bar{\eps}]$, every $y_0\in \Delta_{i^*,\alpha}$, and every $t\geq t_0$. 
We have just checked hypothesis~\ref{hyp:1} of Theorem~\ref{prop-claudia-ext-SW} with $\mathcal{U}=\mathcal{C}$ (remind that $\mathcal{C}=\Delta_{i^*,\alpha}$) and $\mathcal{V}=\R_+^{n}$. Note that the point $y^*$ lies on the boundary of $\R_+^n\times (0,s_{in})$ since $y_i^*=0$ for every $1 \leq i \leq n$ such that $i\not=i^*$ and that 
$y^*_{n+1}=s^*=\vphi(u)\in (0,s_{in})$. But, thanks to the expression defining $f$ as given in \eqref{def-f}, $f$ 
can be straightforwardly extended to a new function (still denoted by $f$) defined over $\R^n \times [0,s_{in}] \times \R_+$. Hence, we are in a position to apply   
Theorem~\ref{prop-claudia-ext-SW} (see \cite[Remark 2.1]{SW99}).  
It follows that there exists $\bar{\eps}_0\in(0,\bar{\eps})$ such that for every $\eps\in [0,\bar{\eps}_0]$, there is a unique point 
$\hat{y}(\eps)\in \Delta_{i^*,\alpha}$ satisfying $f(\hat{y}(\eps),\eps)=0$ and such that $y(t,y_0,\eps)\overset{t\to+\infty}{\longrightarrow} \hat{y}(\eps)$ for every $y_0\in \Delta_{i^*,\alpha}$. We conclude by taking $\eps_{u,\alpha}:=\bar{\eps}_0$ and  $(x^{\eps,u},s^{\eps,u}):=\hat{y}(\eps)$ which, by continuity w.r.t.~$\eps$, is not the wash-out.   
\end{proof}
The previous result asserts a global stability type property for \eqref{sys1} once the dilution rate $u>0$ and the corresponding subset $\Delta_{i^*,\alpha}$ have been chosen where $i^*$ is such that $\mu_{i^*}^{-1}(u)=\vphi(u)<+\infty$ and $\alpha\leq s_{in}-\vphi(u)$. It is also natural to ask how   
\eqref{sys1} behaves if instead we are given some initial condition in $\Delta$. This question can be answered in a simple way as follows. 

\begin{coro} Suppose that the hypotheses of Theorem \ref{main1} are fulfilled and let $y_0:=(x_0,s_0)\in \Delta$ be given such that the 
$i^*$-th coordinate of $x_0$ is positive. Then, 
there is $\eps_{u,y_0}$ such that for every $0\leq \eps \leq \eps_{u,y_0}$, 
the unique solution $y(\cdot)$ to \eqref{sys1} starting at $y_0$ at time $t=0$ converges to $(x^{\eps,u},s^{\eps,u})$. Furthermore, 
every solution $\tilde y$ to \eqref{sys1} such that $\tilde y_{i^*}(0) \geq y_{i^*}(0)$ also converges to $(x^{\eps,u},s^{\eps,u})$. 
\end{coro}
\begin{proof}
Let $\alpha>0$ be such that $\alpha \leq \min(y_{i^*}(0),s_{in}-\vphi(u))$. It follows that $y_0$ belongs to $\Delta_{i^*,\alpha}$ and so is $\tilde y(0)$. 
The result then follows from Theorem \ref{main1}. 
\end{proof}

\subsection{Comments about Theorem~\ref{main1}}
Theorem~\ref{main1} sets up the existence of a perturbed equilibrium point 
$(x^{\eps,u},s^{\eps,u})\not= E_{wo}$  together with its stability in $\Delta_{i^*,\alpha}$.  
It could be interesting to find general conditions on $h$ which guarantee that the perturbed steady-state is such that all species are present asymptotically ({\it{i.e.}},  $x^{\eps,u}_i>0$ for all $i$). 
In the next section, we answer to this question when the perturbation term $h$ is linear w.r.t.~$x$.  

The reasoning as deployed in the proof of Theorem~\ref{main1} does not allow us to obtain global stability of the steady-state in $\Delta\backslash \{0_{\R^n}\}$ for a fixed $\eps>0$. Indeed, for every initial condition $x(0)\in \Delta$ such that  $x_{i^*}(0)=0$, the system \eqref{sys1} with $\eps=0$ 
does not converge to $E_{i^*}$ (but to some other steady state, thanks to the CEP). It follows that  for $\eps=0$, the time to reach 
$E_i^*$ from an initial condition in $\Delta$ such that $x_{i^*}(0)>0$ and $x_{i^*}(0)\rightarrow 0$ goes to infinity. 
This prevents us to apply Theorem~\ref{prop-claudia-MG}  that requires a compact set $\mathcal{C}$ for initial conditions $x(0)$ in order to uniformly control the distance between perturbed and non-perturbed trajectories for all $x(0)\in \mathcal{C}$.  
However, numerical simulations (as performed in Section~\ref{sec:simu}) indicate that global stability of the perturbed system (for a fixed $\eps>0$) 
should be verified in $\Delta\backslash\{0_{\R^n}\}$ in two particular cases for the perturbation term.  Proving global stability in the set $\Delta\backslash\{0_{\R^n}\}$ for a general function $h$ seems rather a delicate question in general and it could be first addressed when the perturbation term is as in Section~\ref{sec:3}.


\section{Analysis of steady-states in the case of a linear exchange term}
\label{sec:3}

The objective of this section is to study properties of steady-states of \eqref{sys1} when the interaction term is linear w.r.t.~$x$ 
and the perturbation parameter enters linearly into the system, that is, 
\begin{equation}{\label{def-h}}
h(x,s,\eps)=\eps T(s)x,
\end{equation}
for all $(x,s,\eps)\in \R^{n+2}$ where $T:\R\rightarrow \R^{n \times n}$ is continuous.  
Considering such a linear coupling w.r.t.~$x$ is motivated by several application models, see, {\it{e.g.}}, \cite{BCLC, deleen,deleen2}.  The fact that the interaction term is now more specific than in the preceding section will allow us to give additional properties concerning 
the steady-state as given in Theorem~\ref{main1}. 
When $h$ is given by \eqref{def-h}, \eqref{sys1} can be  equivalently rewritten as : 
\begin{equation}{\label{sys3}}
\left\{
\begin{array}{cl}
\dot{x}&=B(s,u,\eps)x, \quad  \quad 1 \leq i \leq n,\vspace{0.1cm} \\
\dot{s}&=\ds -\sum_{j=1}^n \frac{\mu_j(s)x_j}{Y_j}+u(s_{in}-s),
\end{array}
\right.
\end{equation}
where the matrix $B(s,u,\eps)\in \R^{n\times n}$ is defined as
$$
B(s,u,\eps):=D(s)-u I_n+\eps T(s).
$$ 
If a pair $(x,s)\in \R^n_+\times [0,s_{in}]$ is a steady state of \eqref{sys3}, then, we necessarily have
\begin{equation}{\label{eq:equi}}
\left\{
\begin{array}{rl}
B(s,u,\eps)x&= 0,  \\
\ds \sum_{j=1}^n \frac{\mu_j(s)x_j}{Y_j}&=u(s_{in} -s). 
\end{array}
\right.
\end{equation}
As a consequence, $E_{wo}$ is always a steady-state of \eqref{sys3}. On the other hand, if a steady-state of \eqref{sys3} is such that $x$ is non-null, then, $x$ is necessarily an eigenvector associated with the $0$ eigenvalue of the matrix $B(s,u,\eps)$. 
To study the occurrence of such a steady-state, in the next section, we will recall classical definitions related to the Perron-Frobenius Theorem and give properties of 
the matrices $T$ and $B$.   
\subsection{Application of the Perron-Frobenius Theorem}

A matrix\footnote{As usual, matrices are named using capital letters and coefficients are represented by lower case letters.
} $A\in \R^{n\times n}$ is said to be  {\it{essentially non-negative}} if $a_{i,j}\geq 0$, for all $i\neq j$ and  {\it{irreducible}} if there are $r>0$ and $k\in \mathbb{N}^*$  
such that all entries of $(A+rI_n)^k$ are positive. Moreover $\lambda(A):=\max\{{\rm Re}(\lambda)\; ; \;  {\rm det}(A-\lambda I_{n})=0  \}$ denotes the largest real part of the eigenvalues of $A$. When dealing with $\lambda(A)$, we shall also refer to the {\it{Perron root}} of $A$.   
Let us now introduce some additional hypotheses about the matrix $T$. 
\label{assm:6}
\begin{enumerate}[label={\rm(A\arabic*)}]
\setcounter{enumi}{4}
    \item\label{as:6}\label{as:7} For all $s\in(0,s_{in}]$, $T(s)$
 is essentially non-negative, irreducible and  
 $\sum_{j=1}^n (T(s)x)_j=0$ for all $x\in \R^n$ ;  
    \item\label{as:8} For all $1 \leq i,j\leq n$ such that $i\not=j$, the function $s\mapsto t_{i,j}(s)$ is non-decreasing and there is $\bar\eps>0$ such that for all $1 \leq i \leq n$, $s\mapsto\mu_i(s)+\eps t_{i,i}(s)$ is also increasing for every $\eps \in [0,\bar{\eps}]$.
\end{enumerate}
For future reference, note that in~\ref{as:8}, $s\mapsto\mu_i(s)+\eps t_{i,i}(s)$ corresponds to the i-th entry of $D(s)+\eps T(s)$. From \ref{as:6} we can prove the following property related to $T(\cdot)$.
\begin{propri}\label{lem:as6} Suppose that \ref{as:6} is satisfied. Then, one has $\lambda(T(s))=0$ for every $s\in[0,s_{in}]$.  
\end{propri}
\begin{proof}
Given $1\leq i \leq n$ and $s\in (0,s_{in}]$, 
one has 
$\sum_{j=1}^n t_{i,j}(s) =
\sum_{j=1}^n \big(T(s)e_i\big)_{j} = 0$, 
where $e_i$ is the ith vector of the canonical basis of $\R^n$. Now, the Perron-Frobenius Theorem implies that
$$ 
0 =  \min_{1 \leq i \leq n} \bigg\{\sum_{j=1}^n t_{i,j}(s) \bigg\} \leq \lambda\big(T(s)^\top\big) \leq \max_{1 \leq i \leq n} \bigg\{\sum_{j=1}^n t_{i,j}(s) \bigg\} =0,
$$
whence the result (using that any matrix and its transpose have  same eigenvalues). 
Since $T$ is continuous, we obtain that $\lambda(T(0))=0$ combining the continuity of 
$\lambda$ (which follows from the continuity of the spectral radius w.r.t.~the matrix) and the fact that $\lambda(T(s))=0$ for any $s\in (0,s_{in}]$. 
\end{proof}

For every $\eps\geq 0$, we define the {\it{critical dilution rate}} as 
\begin{equation}{\label{u-critic}}
u_{\rm c}(\eps):=\lambda\big(D(s_{in})+\eps T(s_{in})\big).
\end{equation}
It will allow us to determine if there is another steady-state than the washout depending on the value of $u$ w.r.t.~$u_c(\eps)$. 
The dependence of the matrix $T$ on $s$ complicates the approach compared to the one presented in \cite{BCLC}. 
Therefore, before delving into this issue, we will first recall two results from \cite{altenberg2012} related to : 

\begin{itemize}
\item[$\bullet$] the convexity of the largest real part of the sum of two matrices\footnote{This result is known as Cohen's convexity Theorem, which pertains to the spectral bound of essentially nonnegative matrices in relation to their diagonal elements.} (Theorem~\ref{altern_thm2} below) ; 
\item[$\bullet$] properties of a function of two variables that is convex in the second variable (Lemma~\ref{altern_lem1} below).  
\end{itemize}
\begin{thm}[\cite{altenberg2012}, Theorem 2]
\label{altern_thm2}
If $D\in\R^{n\times n}$ is a diagonal matrix and $A\in\R^{n\times n}$ an essentially nonnegative matrix, then, $\lambda(A+D)$ is a convex function of $D$.
\end{thm}

\begin{lem}[\cite{altenberg2012}, Lemma 1 (3)]
\label{altern_lem1}
If $f : \R_+^* \times \R_+\rightarrow \R$, $(a,b)\mapsto f(a,b)$ is a function satisfying:
\begin{enumerate}
    \item[{\rm{1}}.] for all $(a,b)\in \R_+^* \times \R_+$ and all $\beta\in \R_+^*$, $f(\beta a,\beta b)=\beta f(a,b)$,
    \item[{\rm{2}}.] for all $a\in \R_+^*$, $b\mapsto f(a,b)$ is convex over $\R_+$,
\end{enumerate}
then, for all $(a,b)\in\R_+^*\times \R_+$, one has 
$\frac{\partial f} {\partial a}(a,b) \leq f(1,0)$ 
except possibly at a countable number of points where the one-sided derivatives\footnote{The notation 
$\frac{\partial f}{\partial a}(a^\pm,b)$ stands for the right and left derivatives of $f$ w.r.t.~$a$ at the point $(a,b)$.} of $f$ exist but differ, and, at these points, one has
$\frac{\partial f}{\partial a}(a^-,b) < \frac{\partial f}{\partial a}(a^+,b)\leq f(1,0)$. 
\end{lem}

The next statement establishes properties of the matrix $B$ and of the critical dilution rate $u_c$ 
extending \cite[Lemma 3.1]{BCLC} to the case of a non-constant and non-necessarily symmetric matrix $T(\cdot)$, 
 and whenever yield coefficients are non-necessarily equal. It is also useful for proving Proposition~\ref{prop:terence}. 
\begin{prop} {\label{lem-sec3}}
Suppose that \ref{as:1} and\ref{as:7}-\ref{as:8} hold true. \smallskip \\
{\rm{(i)}}. For $(\eps,u)\in[0,\bar\eps]\times\R_+^*$, one has $\lambda(B(0,u,\eps))<0$ and $s\mapsto \lambda(B(s,u,\eps))$ is increasing over $[0,s_{in}]$. 
\smallskip \\
{\rm{(ii)}}. The function $u_c(\cdot)$ is non-increasing over $\R_+$. Moreover, one has 
$u_{\rm c}(0)=\mu_{max}$ and $u_{\rm c}(\eps)\overset{\eps\to+\infty}{\longrightarrow}\hat{\mu}$ 
where $\mu_{max}:=\ds \max_{i=1,\ldots,n} \mu_i(s_{in})$, $\hat{\mu}:=w^\top D(s_{in})v$, and $v, w\in \R^n$  are uniquely defined by
\begin{equation}{\label{def-vw}}
T(s_{in})v =0 \; ; \; w^\top T(s_{in})=0 \; ; \; w^\top v=1.  
\end{equation}
\end{prop}

\begin{proof} To prove (i), let $\eps\in [0,\bar \eps]$ and $u>0$. 
One has $B(0,u,\eps)=-uI_n+\eps T(0)$, thus, 
$\lambda(B(0,u,\eps))=-u+\eps\lambda(T(0))=-u<0$ (thanks to \ref{as:7}). From the Perron-Frobenius Theorem, for every $s\in (0,s_{in}]$, 
$\lambda(B(s,u,\eps))$ exists and is of multiplicity one, hence, we can uniquely define $v(s)$ as the unitary associated eigenvector (which is thus with positive coordinates). Now, given $0<s<s'$, \ref{as:8} implies that $B(s,u,\eps)\leq B(s',u,\eps)$ entry-wise. 
Thus, by the Collatz-Wielandt formula\footnote{For every essentially non-negative matrix $A$, $\lambda(A)= \sup_{x\in\R^n_+\backslash\{0\}}\min_{1 \leq i \leq n} \frac{(Ax)_i}{x_i}$.} and using \ref{as:8}, we deduce that   
\begin{align*}
    \lambda(B(s',u,\eps)) &= \sup_{x\in\R^n_+\backslash\{0\}}\min_{1 \leq i \leq n} \frac{(B(s',u,\eps)x)_i}{x_i} \\ &\geq \min_{1 \leq i \leq n} \frac{(B(s',u,\eps)v(s))_i}{v_i(s)} 
    >\min_{1 \leq i \leq n} \frac{(B(s,u,\eps)v(s))_i}{v_i(s)}  = \lambda(B(s,u,\eps)) , 
\end{align*}
where the strict inequality comes from the second part of \ref{as:8}. We obtain that $s\mapsto \lambda(B(s,u,\eps))$ is strictly monotone over $(0,s_{in}]$ and, therefore, over $[0,s_{in}]$ as well which proves (i). 

To prove (ii), we adapt the proof of Proposition 3.5 in \cite{BCLC}  
using the aforementioned results (\cite{altenberg2012}) in the case where $T$ depends on $s$. 
First, we can directly conclude from the CEP that $u_{\rm c}(0)=\mu_{\max}$. Recall now that for any essentially non-negative and irreducible matrix $A$, the Perron-Frobenius Theorem ensures that $\lambda(A)$ is simple, hence, there is a neighborhood $\mathcal{V}_A$ of $A$ in $\R^{n\times n}$ such that $\lambda(B)$ is simple for every $B\in \mathcal{V}_A$.  
Furthermore, the following (standard) result holds true (see, {\it{e.g.}}, \cite{DN84}) : 
the derivative of $\lambda$ w.r.t.~$A$, denoted by $D_1^A$ is given by $D_1^A := w_a v_a^\top$,
where $v_a,w_a\in\R^n$ are uniquely defined by   
$$(\lambda(A)I_n-A)v_a =0, \qquad w_a^\top (\lambda(A)I_n-A)=0 \quad \text{and}\quad w_a^\top v_a=1. $$
Observe that $T(s_{in})$ is essentially non-negative, irreducible and that $\lambda(T(s_{in}))=0$. 
Thus, if $v$ and $w$ are given by \eqref{def-vw}, one has $D_1^{T(s_{in})} = wv^\top$. 
Next, let us set 
$$
f(a,b):=\lambda(a T(s_{in})+b D(s_{in}))
$$
for $(a,b)\in \R_+ \times \R$ and let $\gamma(b):=f(1,b)$.  
When $b \downarrow 0$, the Perron root $\lambda(T(s_{in})+bD(s_{in}))$ is simple since $\lambda(T(s_{in}))$ is also simple, hence, $\gamma$ is differentiable at $b=0$ and using the chain rule formula, we get that  
$$
  \gamma'(0)= \sum_{1 \leq i, j \leq n} (w v^\top)_{i,j} \ (D(s_{in}))_{i,j} = w^\top D(s_{in}) v.
$$
Thus, using a first-order expansion of $\gamma$ around $0$, we obtain as $\eps\rightarrow +\infty$: 
$$ 
u_{\rm c}(\eps)= \eps \lambda\left(T(s_{in})+ \frac{1}{\eps}D(s_{in})\right) = \eps \left( \lambda (T(s_{in}))+\frac{1}{\eps} w ^\top D(s_{in}) v + 
o\left(\tfrac{1}{\eps}
\right)\right) = w ^\top D(s_{in}) v+ o(1), 
$$
where $o(1)\rightarrow 0$ as $\eps\rightarrow +\infty$.  
We deduce that $u_{\rm c}(\eps)\to \hat{\mu}$ as $\eps\rightarrow +\infty$. To show that $u_{\rm c}$ is non-increasing over $\R_+$,  
we use Theorem~\ref{altern_thm2} and Lemma~\ref{altern_lem1}. Let us first check that $f$ fulfills the hypotheses of Lemma \ref{altern_lem1}. 
Doing so, observe that for $\beta>0$, $f$ satisfies 
$f(\beta a,\beta b) =  \beta f(a,b)$. 
Next, Theorem \ref{altern_thm2} implies that $D\mapsto \lambda(T(s_{in})+D(s_{in}))$ is convex w.r.t.~the diagonal matrix $D$ given that by Assumption \ref{as:6}, $T(s_{in})$ is essentially non-negative. This property implies that $f$ is convex w.r.t.~$b$, hence, we are in a position to apply Lemma~\ref{altern_lem1}. 

Furthermore, from \ref{as:7}, the matrix $D(s_{in})+\eps T(s_{in})$ is essentially non-negative and irreducible for all $\eps>0$ so that $\lambda(D(s_{in})+\eps T(s_{in}))$ is simple. It follows that 
$\eps \mapsto u_c(\eps)$ is differentiable over $\R_+^*$ and that for every $\eps >0$, one has:
$$
u'_c(\eps)  = \frac{\partial f}{\partial \eps}(\eps,1) \leq f(1,0) =\lambda(T(s_{in}))  = 0,
$$
and therefore, $u_c$ is non-increasing\footnote{We are ignoring that the one-side derivatives may differ (see Lemma \ref{altern_lem1}) since $u_{\rm c}$ is differentiable.}.

\end{proof}
As a consequence, under \ref{as:1}, \ref{as:7}, and \ref{as:8}, $\eps \mapsto u_{\rm c}(\eps)$ is one-to-one from $\R_+$ into $(\hat{\mu},\mu_{\max}]$, hence its reciprocal $\xi:(\hat \mu,\mu_{max}]\rightarrow \R_+$ is defined as
$$
\xi(u):=
\left\{
\begin{array}{lcl}
+\infty & \mathrm{if} & u \in [0,\hat{\mu}],\\
u_{\rm c}^{-1}(u) & \mathrm{if} & u \in (\hat{\mu},\mu_{\max}].
\end{array}
\right.
$$
\subsection{Coexistence steady-state}
In this section, we discuss the occurrence of a coexistence steady-state when the perturbation term is given by \ref{def-h}.  
Our discussion on steady-states of \eqref{sys3} is divided into two cases, whether $u\geq u_{\rm c}(\eps)$ or $u<u_{\rm c}(\eps)$. The next proposition presents similarities with \cite[Propostion 3.3]{BCLC}, but, as explained before, 
yield coefficients are not necessarily all equal to one (in contrast with \cite{BCLC}) and the perturbation matrix $T(s)$ is now a function of $s$ (and non-necessarily symmetric). That is why we provide the proof in details.
\begin{prop}\label{prop:terence}
Suppose that \ref{as:1} and \ref{as:7}-\ref{as:8} hold true and that there is a unique $i^*\in\{1,...,n\}$ such that $\mu^{-1}_{i^*}(u)=\vphi(u)<+\infty$. 
\smallskip
\\
{\rm(i)} If $(\eps,u)\in [0,\bar \eps]\times \R_+^*$ is such that $u\geq u_{\rm c}(\eps)$, then $E_{wo}$ is the only equilibrium of \eqref{sys3} and it is stable. If $u> u_{\rm c}(\eps)$, then $E_{wo}$ is globally asymptotically stable. 
\smallskip
\\
{\rm(ii)} If $(\eps,u)\in [0,\bar \eps]\times \R_+^*$ is such that $u< u_{\rm c}(\eps)$, then, there are only two equilibria for \eqref{sys3}, namely the wash-out $E_{wo}$ that is unstable and a coexistence steady-state $E_{\eps,u}:=(x^{\eps,u},s^{\eps,u})\in (0,+\infty)^n \times (0,s_{in})$. Furthermore, for every $u\in (0,\mu_{\max})$, there is $\tilde{\eps}_u\in (0, \tilde \xi(u))$ (where $\tilde \xi(u):=\min(\bar \eps,\xi(u))$) such that for every $\eps\in (0,\tilde \eps_u)$, $E_{\eps,u}$ is locally asymptotically stable.  
\end{prop}

\begin{proof}
Assume that $u\geq u_{\rm c}(\eps)$. We will see that $E_{wo}$ is the only solution of \eqref{eq:equi}. Suppose by contradiction that there is 
$(\bar x, \bar s)\in \R_+^n \times [0,s_{in}]$ satisfying \eqref{eq:equi} and such that $\bar x\not= 0_{\R^n}$. Hence, we necessarily have $\bar s<s_{in}$. Moreover, $\bar x$ is an eigenvector of $B(\bar s,u,\eps)$ associated with the eigenvalue $0$, implying that $\lambda(B(\bar s,u,\eps))\geq 0$. On the other hand, by Proposition~\ref{lem-sec3}, $s\mapsto\lambda(B(s,u,\eps))$ is increasing, thus 
$$
\lambda(B(\bar s,u,\eps)) < \lambda(B(s_{in},u,\eps)) = \lambda(D(s_{in})+\eps T(s_{in}))-u \leq \lambda(D(s_{in})+\eps T(s_{in}))-u_{\rm c}(\eps) = 0,
$$
therefore we have a contradiction so that we must have $(\bar x,\bar s)=E_{wo}$. 
Let us now turn to the stability of $E_{wo}$. We claim that for every $\eta>0$ there is a $\delta>0$ such that $|{x(t)}|<\eta$ and $|s(t)-s_{in}|<\eta$ for every solution $(x(t),s(t))$ of \eqref{sys3} with $\norm{x(0)}\leq \delta$ and $|s_{in}-s(0)|\leq \delta$. First, notice that $\dot{x}=(D(s)-u I_n+\eps T(s))x\leq (D(s_{in})-u I_n+\eps T(s_{in}))x$, but $\lambda(B(s_{in},u,\eps))\leq 0$, hence, we deduce that $|x(t)| \leq  |x(0)|$ for every $t\geq 0$. Now, fix $\eta>0$ and observe that 
$$
\dot{s}(t)\geq - \left(\sum_{j=1}^n \frac{\mu_j(s(t))^2}{Y_j^2}\right)^{1/2}|x(t)|
+u(s_{in}-s(t))\geq -C |x(0)|+u(s_{in}-s(t)), \quad t \geq 0,  
$$
where $C:=\big{(}\sum\frac{\mu_j(s_{in})^2}{Y_j^2}\big{)}^{1/2}$. Let us set $\delta:=\min(\eta, u\eta/(2C))$ and then take $|x(0)|\leq \delta$ and $|s_{in}-s(0)|\leq \delta$. We deduce first that $|x(t)|\leq \delta \leq \eta$ for all $t\geq 0$ and that if at some time $t'\geq 0$ we have $s(t')=s_{in}-\eta$, then
$$
\dot{s}(t') \geq -\frac{uC\eta}{2C}+u \eta=\frac{u \eta}{2}>0. 
$$
It follows that for all time $t\geq 0$, one has $s(t) \geq s_{in}-\eta$ which allows us to conclude on the stability of $E_{wo}$.  
Finally, in the case where $u> u_c(\eps)$, we have
$$
\frac{d}{dt}|x(t)|^2\leq 2 x(t)^\top B(s(t),u,\eps)x(t) \leq 2 x(t)^\top B(s_{in},u,\eps)x(t),
$$
using that $\eps\in [0,\bar \eps]$. Since $\lambda(B(s_{in},u,\eps))<0$, we conclude that $|x(t)| \rightarrow 0$ as $t\rightarrow +\infty$ which proves (i). 

Suppose now that $0<u<u_{\rm c}(\eps)$ and let $(x,s)$ be a steady-state of \eqref{sys3}. If $ s=s_{in}$, then, we have $ x=0$. We easily check that the Jacobian matrix of \eqref{sys3} at $E_{wo}$ is a block matrix whose first block is $B(s_{in},u,\eps)$ which is such that $\lambda(B(s_{in},u,\eps))=u_{\rm c}(\eps)-u>0$. It follows that $E_{wo}$ is unstable, thus, $s\in (0,s_{in})$.  
Thanks to the Perron-Frobenius Theorem, $x$ is an eigenvector associated with the zero eigenvalue of $B(s,u,\eps)$ (the largest one), hence, it is with positive coordinates.  

The vector $x$ is uniquely defined as follows. Since $s\mapsto \lambda(B(s,u,\eps))$ is continuous and increasing and given that $\lambda(B(0,u,\eps))=-u<0$ and $\lambda(B(s_{in},u,\eps))>0$, there is a unique solution $s:=s^{\eps,u}\in (0,s_{in})$ to the equation $\lambda(B(s,u,\eps))=0$. 
Now, $x$ is proportional to the Perron vector\footnote{It is the unique unitary eigenvector in the eigenspace of $B(s,u,\eps)$ for the zero eigenvalue.} $a^{\eps,u}$  associated with the zero eigenvalue of $B(s,u,\eps)$, thus, there is a unique $\alpha^{\eps,u} \in \R$ such that $ x=\alpha^{\eps,u} a^{\eps,u}$. From \eqref{eq:equi}, we conclude that 
\begin{equation}{\label{def-scalar}}
\alpha^{\eps,u}:=\frac{u(s_{in}-s^{\eps,u})}{\tau} \; \; \mathrm{where} \; \; \tau:=\sum_{j=1}^n \frac{\mu_j(s^{\eps,u})a_j^{\eps,u}}{Y_j}.  
\end{equation}
The vector $x^{\eps,u}:=x$ is then uniquely defined. Finally, it is well-known that for $\eps=0$, the Jacobian matrix of \eqref{sys3} at the steady-state is a Hurwitz matrix (see, {\it{e.g.}}, \cite{alain-livre,SmithWalt}). The existence of $\tilde \eps_u\leq \tilde \xi(u)$ follows from the continuity of the eigenvalues of a matrix w.r.t.~parameters. This ends the proof. 
\end{proof}
\begin{rem}{\label{rem-3.2}}
{\rm{(i)}}.  The equilibrium $E_{\eps,u}=(x^{\eps,u},s^{\eps,u})$ is called {\it{coexistence steady-state}} because one has $x^{\eps,u}_i>0$ for every $1 \leq i \leq n$. Remind that $s^{\eps,u}$ is defined by $\lambda(B(s,u,\eps))=0$ and that $x^{\eps,u}=\alpha^{\eps,u}a^{\eps,u}$ where $\alpha^{\eps,u}$ is defined in \eqref{def-scalar} and $a^{\eps,u}$ is the (unitary) Perron vector associated with $\lambda(B(s,u,\eps))$. Note that for an arbitrary pair $(\eps,u)\in \R_+^* \times [0,u_c(\eps))$, it would be interesting to know whether the Jacobian matrix of \eqref{sys2} at this point is Hurwitz, however, this question presents some difficulties. 
\smallskip
\\
{\rm{(ii)}}. In propositions~\ref{lem-sec3} and \ref{prop:terence}, the fact that $\eps\in [0,\bar \eps]$ is not optimal (in the sense that we only conclude 
for $\eps\in [0,\bar \eps]$), but, the parameter $\bar \eps$ is not necessarily 
``small'' since it is defined by $T$ (for instance, $\bar \eps=1/2$ in the example of section~\ref{case2}). 
When $T$ does not depend on $s$, \ref{as:8} is obviously true for all $\eps\geq 0$, hence, $\bar \eps=+\infty$. 
Finally, extending these results for $\eps> \bar \eps$ when $T$ depends on $s$ is certainly not straightforward, and it may not even be true in general.  
\end{rem}


\section{Numerical simulations}
\label{sec:simu}

In this section, we provide a numerical analysis of \eqref{sys1}, our goal being to confirm the theoretical study conducted previously while highlighting other properties of the system numerically. This numerical study is conducted for two types of linear perturbations w.r.t.~$x$: 
\begin{itemize}
\item[$\bullet$] Case 1 : $h(x,s,\eps):= \eps \Theta x$ where $\Theta\in \R^{n\times n}$ is a given symmetric matrix ; 
\item[$\bullet$] Case 2 : $h(x,s,\eps):= \eps T(s) x$ where for every $s\in [0,s_{in}]$, $T(s)\in \R^{n\times n}$ is a given matrix involving the kinetics $\mu_i$ and that is non-symmetric neither constant. 
\end{itemize} 
For case 1, we shall refer to the case with a constant perturbation matrix (or transition matrix) and for case 2, we shall refer to a non-constant perturbation matrix (or transition matrix). The exact definition of $\Theta$ and $T(s)$ and the discussion of the corresponding dynamical system can be found in Section~\ref{case1} and~\ref{case2} respectively.  
Section~\ref{subsec:num} presents numerical simulations of trajectories that are carried out with the SciPy library in \texttt{Python} and more specifically using the function \href{https://docs.scipy.org/doc/scipy/reference/generated/scipy.integrate.solve_ivp.html}{\texttt{solve\_ivp}} that solves numerically ordinary differential equations using an explicit Runge-Kutta method of fifth order\footnote{The scripts for reproducing the numerical experiments of this paper can be found in the repository 
\url{https://github.com/calvarezlatuz/simulation_CM}}. The number of grid points is around 400. 

\begin{figure}[H]
    \centering
    \includegraphics[width=0.5\linewidth]{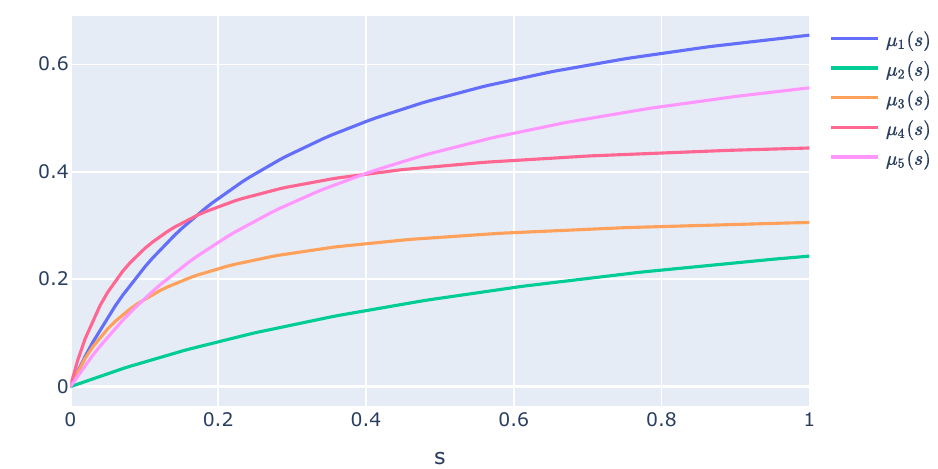}
    \caption{Plot of the kinetics $\mu_i(s)= \frac{a_i s}{b_i+s}$ of Monod type  for $i=1,...,5$ (coefficients are given in Table \ref{tab:parameters}).}
    \label{fig:Monod_kinetics}
\end{figure}

The parameters used for the numerical simulations are defined as follows.  
We assume that kinetics are of Monod type, {\it{i.e.}}, 
$\mu_i(s)= \frac{a_i s}{b_i+s}$,  
where $a_i, b_i$ are positive parameters for every $1 \leq i \leq n$ (it is easily seen that  \ref{as:1} is fulfilled). 
Figure \ref{fig:Monod_kinetics} depicts the kinetics that have been chosen for the examples. The values of yield coefficients and other coefficients can be found below. 

In the numerical simulations of \eqref{sys1} related to its asymptotic behavior, note that we have chosen a constant dilution rate value fixed to $u=0.4$. 
From the CEP, we deduce that for $\eps=0$, then species $1$ survives (see Fig.~\ref{fig:Monod_kinetics}) so that in Theorem~\ref{main1}, 
the index $i^*$ is such that $i^*=1$. 
{\renewcommand{\arraystretch}{1.5}
\begin{table}[H]
    \centering
    \begin{tabular}{llccccc}
        \toprule
        \multicolumn{1}{l}{\textbf{Species number }} &$i$& \textbf{1} & \textbf{2} & \textbf{3} & \textbf{4} & \textbf{5}  \\ 
        \midrule
        \multirow{2}{*}{\textbf{Coefficients of $\mu_i$}} & $a_i$ & 0.84 & 0.46 & 0.34 & 0.48 & 0.76 \\ 
                                       & $b_i$ & 0.28 & 0.90 & 0.11 & 0.09 & 0.36 \\ 
        \midrule
        \multirow{1}{*}{\textbf{Yield coefficients}} & $Y_i$ & 1.00 & 1.50 & 2.00 & 2.50 & 3.00 \\ 
        \midrule\midrule 
        \multicolumn{2}{l}{\textbf{Simulation parameters}} 
                                       & $n=5$        & $s_{in}=1$    & $u=0.4$ & $\eps=0.1$ &  $\alpha=0.05$  \\ 
        \bottomrule
    \end{tabular}
    \caption{Values of $a_i$, $b_i$, $Y_i$, and simulation parameters.}
    \label{tab:parameters}
\end{table}}

\subsection{Case 1 : perturbation term defined via a constant perturbation matrix}{\label{case1}}
In this section, we suppose that $h(x,s,\eps):= \eps \Theta x$ where the matrix $\Theta$ is defined as
{\small{$$ \Theta := \left[\begin{array}{ccccc} -1 & 1 & 0 & \cdots & 0\\ 1 & -2 & 1 & \cdots & 0\\ \vdots & \ddots & \ddots & \ddots & \vdots\\ 0 & \cdots & 1 &  -2 & 1\\ 0& \cdots & 0& 1& -1 \end{array}\right]$$}}

so that \eqref{sys1} rewrites
\begin{equation}
\label{eq:lineal-CS}
\left\{
    \begin{array}{rlll}
        \dot{x} &= D(s)x-u x+\varepsilon \Theta x, \\
        \dot{s} &= - \displaystyle\sum_{j=1}^n\mu_j(s)\frac{x_j}{Y_j}+ u(s_{in}-s).  
    \end{array}
\right.
\end{equation}

To explain, the origin of $\Theta$, let us examine the following population dynamics model structured around a single phenotypic trait called $z$ and 
where $z$ is with values within a fixed segment $I\subset \R$. 
It describes the evolution of both species and substrate with concentrations $x(t,z)$ and $s(t,z)$ respectively: 
\begin{equation}
\label{eq:lineal-CS-continuity}
\left\{
    \begin{array}{rlll}
        \partial_t{x} (t,z)&= \mu(s(t),z)x(t)-u x(t)+\Delta_z x(t,z),\smallskip \\
        \dot{s} (t)&= - \int_I \frac{\mu(s(t),z')}{Y(z')}+ u(s_{in}-s(t)), 
    \end{array}
\right.
\end{equation}
where $t>0$ and $z\in I$. 
In the preceding integro-differential system, $\mu$ denotes the kinetics (depending now on the substrate and on the trait), $\Delta_z x$ denotes the Laplacian of $x$ w.r.t.~$z$, and the function $z\mapsto Y(z)$ replaces the yield coefficients $Y_i$ to be found in the finite-dimensional setting for the species.  
To ensure that the integro-differential system \eqref{eq:lineal-CS} is well-posed, it is typically coupled with appropriate initial and boundary conditions, which help to guarantee both the uniqueness and existence of a solution to \eqref{eq:lineal-CS}. 
Since our purpose in this paper is to address stability properties for a finite number of species, we will not further discuss the preceding system. 
We just want to point out that $\Theta$ can be seen as a discretization of the diffusion term $\Delta_z x$ for Neumann boundary conditions.

Note that the main difference between \eqref{eq:lineal-CS} and the model studied in \cite{BCLC} is that in the present setting, yield coefficients are not all equal to one (in contrast with the model in \cite{BCLC}). Furthermore, it is important to note that the dynamical system examined in \cite{BCLC} is only analyzed in the asymptotic case where $u$ approaches zero. In contrast, our study maintains $u$ at a fixed constant value while $\eps$ is a small parameter (which seems more reasonable from a practical viewpoint). 

In order to apply the results of Section~\ref{sec:3}, we can easily 
verify the following : 
\begin{itemize}
\item[$\bullet$] the matrix $\Theta$ is symmetric, essentially non-negative, irreducible and such that 
$\sum_{j=1}^n (\Theta x)_j=0$ for every $x\in \R^n$ (see also \cite{BCLC}) so that \ref{as:6} is satisfied.   
\item[$\bullet$] assumption \ref{as:8} is satisfied for every $\eps\geq 0$ so that $\bar \eps=+\infty$. 
\end{itemize}
Hence, for every $(\eps,u)\in \R_+^*\times \R_+^*$ such that $u<u_c(\eps)$, Section~\ref{sec:3} ensures the existence of a unique coexistence steady-state 
$(x^{\eps,u},s^{\eps,u})$ of \eqref{eq:lineal-CS} apart from the wash-out. Global stability of this steady-state is described in Theorem~\ref{main1}. 

Finally, it can be checked that if $x(0)\not=0$, then, the corresponding solution to \eqref{eq:lineal-CS} satisfies $x_i(t)>0$ for every $t>0$ and every $1 \leq i \leq n$ (to prove this, it is sufficient to apply Lemma~\ref{lem-positivity}, whose main hypothesis \eqref{CS-pos} is easy to verify). 
In view of the link between $\Theta$ and the discretization of the Laplacian operator, this property is in line with the well-known regularizing effect of the heat equation.


\subsection{Case 2 : perturbation term defined via a non-constant perturbation matrix}{\label{case2}}
In this section, we suppose that $h(x,s,\eps):= \eps T(s) x$ where the matrix $T(s)$ is defined as
{\small{
\begin{equation}{\label{NC-MATRIX}}
\begin{bmatrix}
-2\mu_1(s) & \mu_2(s) & 0 & \cdots & \mu_n(s) \\
\mu_1(s) & -2\mu_2(s) & \mu_3(s) & \cdots & 0 \\
0 & \mu_2(s) & -2\mu_3(s) & \cdots & 0 \\
\vdots & \vdots & \vdots & \ddots & \mu_{n}(s) \\
\mu_1(s) & 0 & 0 & \mu_{n-1}(s) & -2\mu_{n}(s)
\end{bmatrix}
\end{equation}
}}

so that \eqref{sys1} rewrites
\begin{equation}
\label{eq:cultivated-CS}
\left\{
    \begin{array}{rlll}
        \dot{x} &= D(s)x-ux+\eps T(s)x\\
        \dot{s} &= - \displaystyle\sum_{j=1}^n\mu_j(s)\frac{x_j}{Y_j}+ u(s_{in}-s),
    \end{array}
\right.    
\end{equation}
Observe that $T(s)$ is a tridiagonal matrix that includes two coefficients located at the position $(1,n)$ and $(n,1)$.  
The origin of this matrix is derived from \cite{courgeau1969mutations}, which describes the effect of mutations within a cultivated biomass, leading to the consideration of each species as a distinct gene mutation. For every $s$, the so-called {\it{mutation matrix }} $T(s)$ 
is derived considering a probability of mutation, leading that way to the above circulant matrix. The system as described in \eqref{eq:cultivated-CS} is also closely related to the one studied in \cite[Section~4.2]{lobry}. 

In order to apply the results of Section~\ref{sec:3}, we check \ref{as:6} and \ref{as:8}:
\begin{itemize}
\item[$\bullet$]
it is easily seen that for any $s\in(0,s_{in}]$, the matrix $T(s)$ is essentially non-negative, irreducible and such that $\sum_{j=1}^n (T(s) x)_j=0$ for every $x\in \R^n$ so that \ref{as:6} is satisfied.  
\item[$\bullet$] \ref{as:8} is satisfied for every $\eps\in[0,\bar \eps)$ with $\bar \eps:=1/2$.  
\end{itemize}
Hence, for every $(\eps,u)\in [0,\bar \eps]\times \R_+^*$ such that $u<u_c(\eps)$, Section~\ref{sec:3} ensures the existence of a unique coexistence steady-state 
$(x^{\eps,u},s^{\eps,u})$ of \eqref{eq:cultivated-CS} apart from the wash-out. As in the previous case, global stability of this steady-state is described in Theorem~\ref{main1} and all species are with positive values over $(0,+\infty]$ provided that $x(0)\not=0$ (the hypothesis \eqref{CS-pos} in Lemma~\ref{lem-positivity} is also easy to check).  

Interestingly, with the data given in Table~\ref{tab:parameters}, we verified numerically that $E_{\eps,u}$ still exists for every 
$\eps>\bar \eps$ and every $u<u_c(\eps)$ (see the discussion in Remark~\ref{rem-3.2} (ii)). 
This is why, the numerical analysis as performed in Section~\ref{subsec:num} 
has also been carried out for $\eps>\bar \eps$ (this is only for~figures~\ref{fig:d0_cult} and~\ref{fig:Hurwitz_cult}). 
The next two figures show (numerically) that $s\mapsto \lambda(B(s,u,\eps)$ is always increasing which is sufficient 
to ensure the validity of Proposition~\ref{lem-sec3} and \ref{prop:terence} when $T$ is given by \eqref{NC-MATRIX}. 
\begin{figure}[H]
    \centering
    \includegraphics[width=0.4\linewidth]{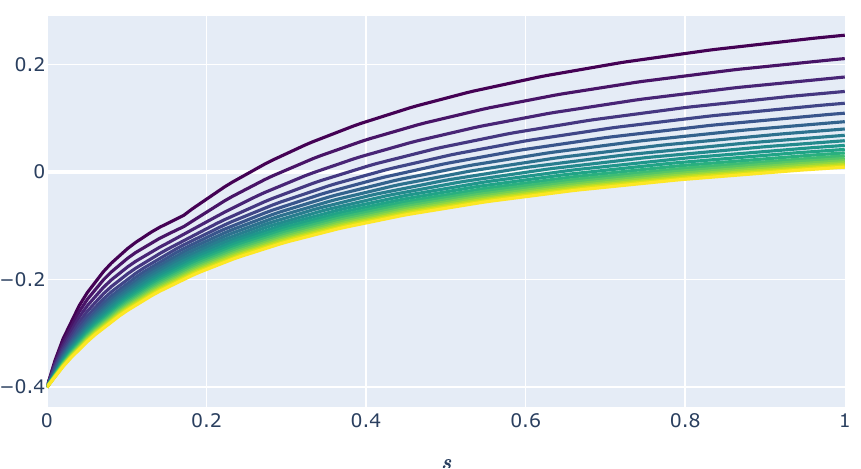}
    \includegraphics[width=0.4\linewidth]{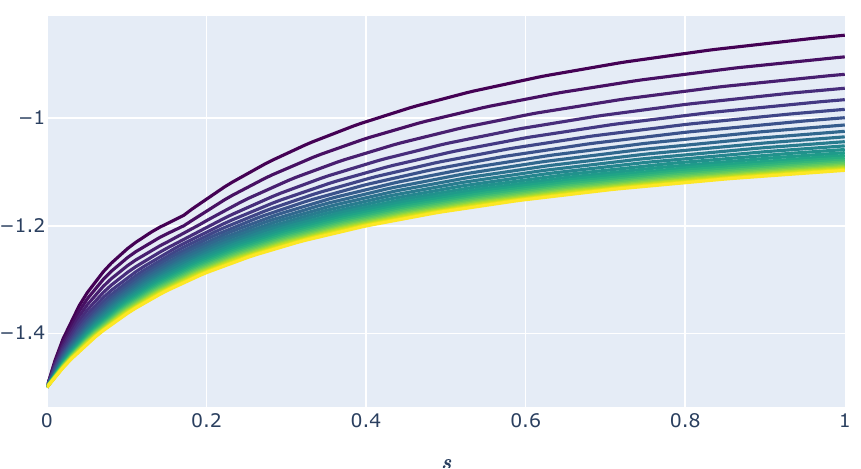}
    \caption{Plot of $s\mapsto \lambda(B(s,u,\eps))$ over $[0,s_{in}]$ for $\eps$ discretized between $0$ and $100$ and for $u=0.4$ (left) and $u=1.5$ (right) in the case where $T$ is given by~\eqref{NC-MATRIX}.}
    \label{fig:verif}
\end{figure}
 

\subsection{Results for both examples}{\label{subsec:num}}

\begin{figure}[H]
    \centering
    \includegraphics[width=0.93\linewidth]{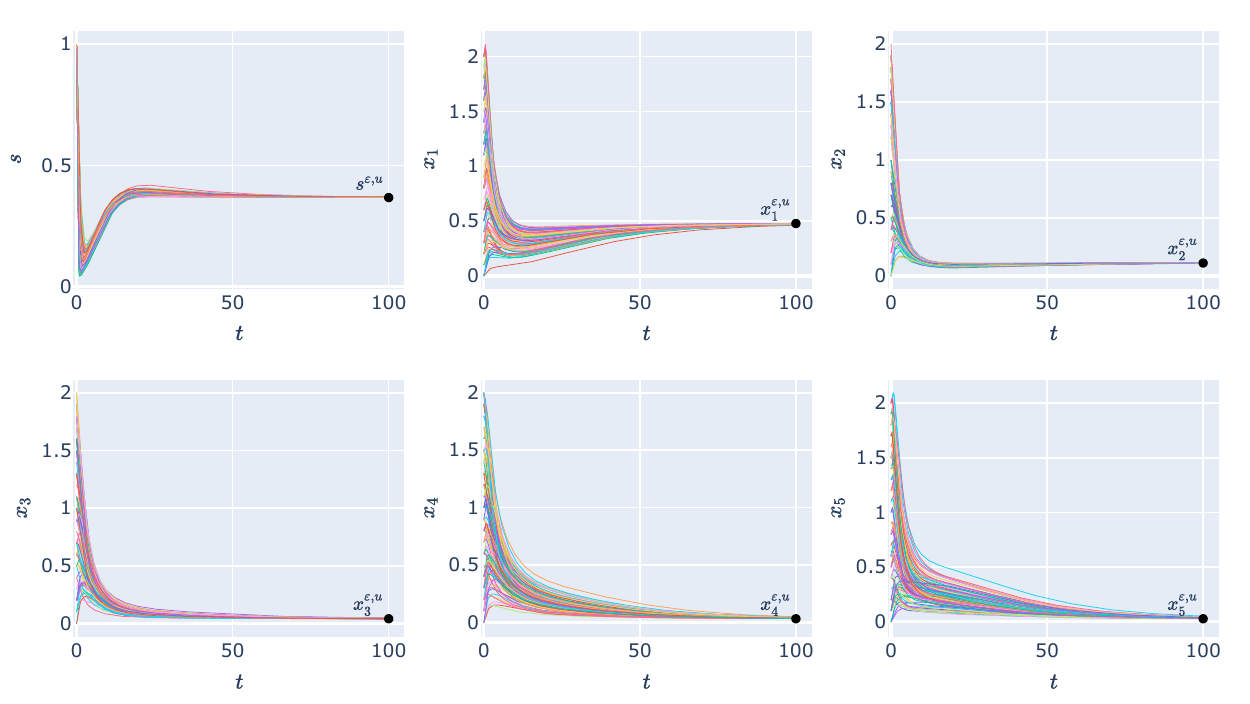}
    \caption{Solutions to \eqref{eq:lineal-CS} for 100 randomly generated initial conditions in the set $(0,10]^n \times [0,s_{in}]$}
    \label{fig:monod_lin1}
\end{figure}

\begin{figure}[H]
    \centering
    \includegraphics[width=0.93\linewidth]{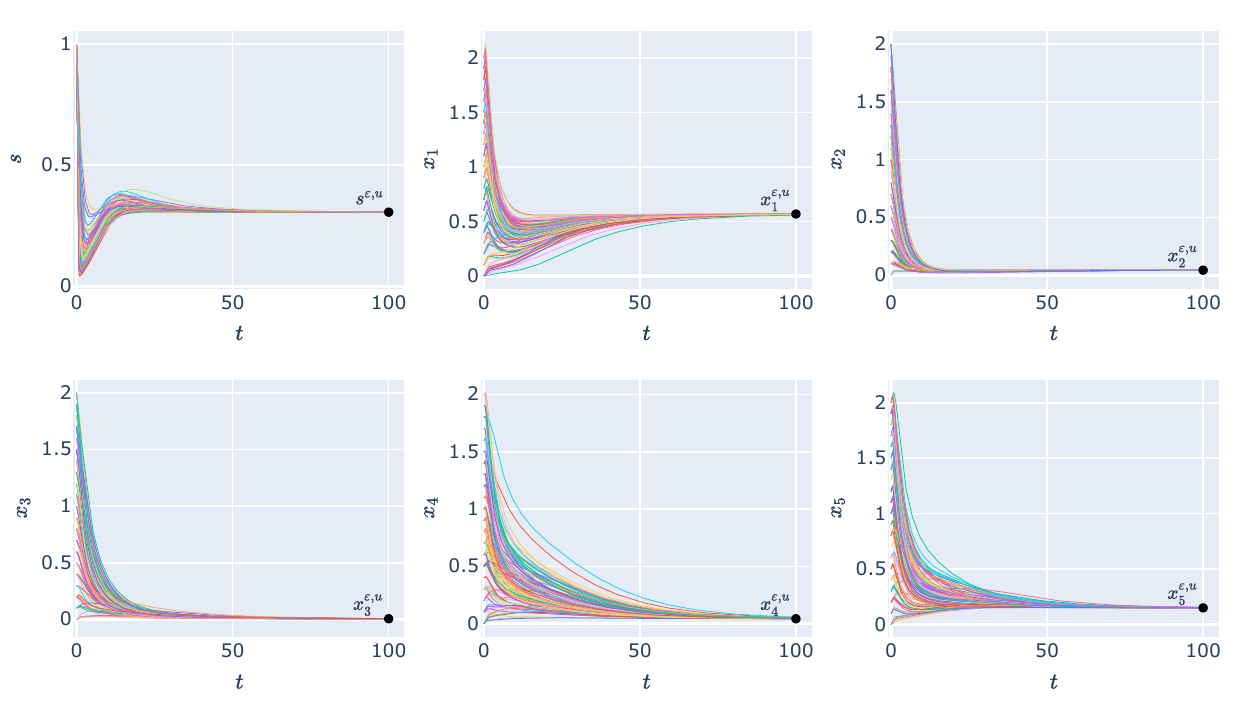}
    \caption{Solutions to \eqref{eq:cultivated-CS} for 100 randomly generated initial conditions in the set $(0,10]^n \times [0,s_{in}]$}
    \label{fig:monod_cul1}
\end{figure}

\begin{figure}[H]
    \centering
    \includegraphics[width=0.99\linewidth]{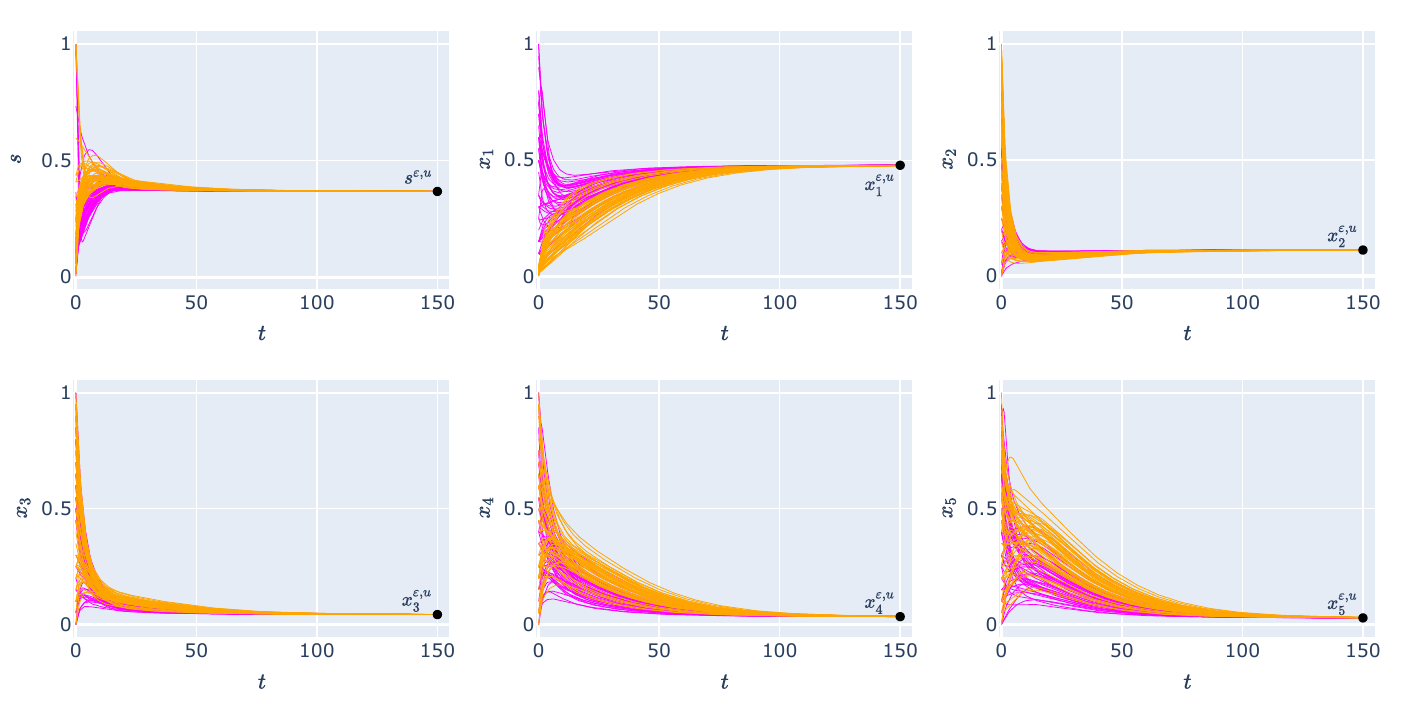}
    \caption{
Solutions to \eqref{eq:lineal-CS} for 50 randomly generated initial conditions in the set $\Delta_{1,\alpha}$ (magenta) and 50 others in the set 
$\Delta\backslash \Delta_{1,\alpha}$ (orange). }
    \label{fig:monod_lin2}
\end{figure}

\begin{figure}[H]
    \centering
    \includegraphics[width=0.99\linewidth]{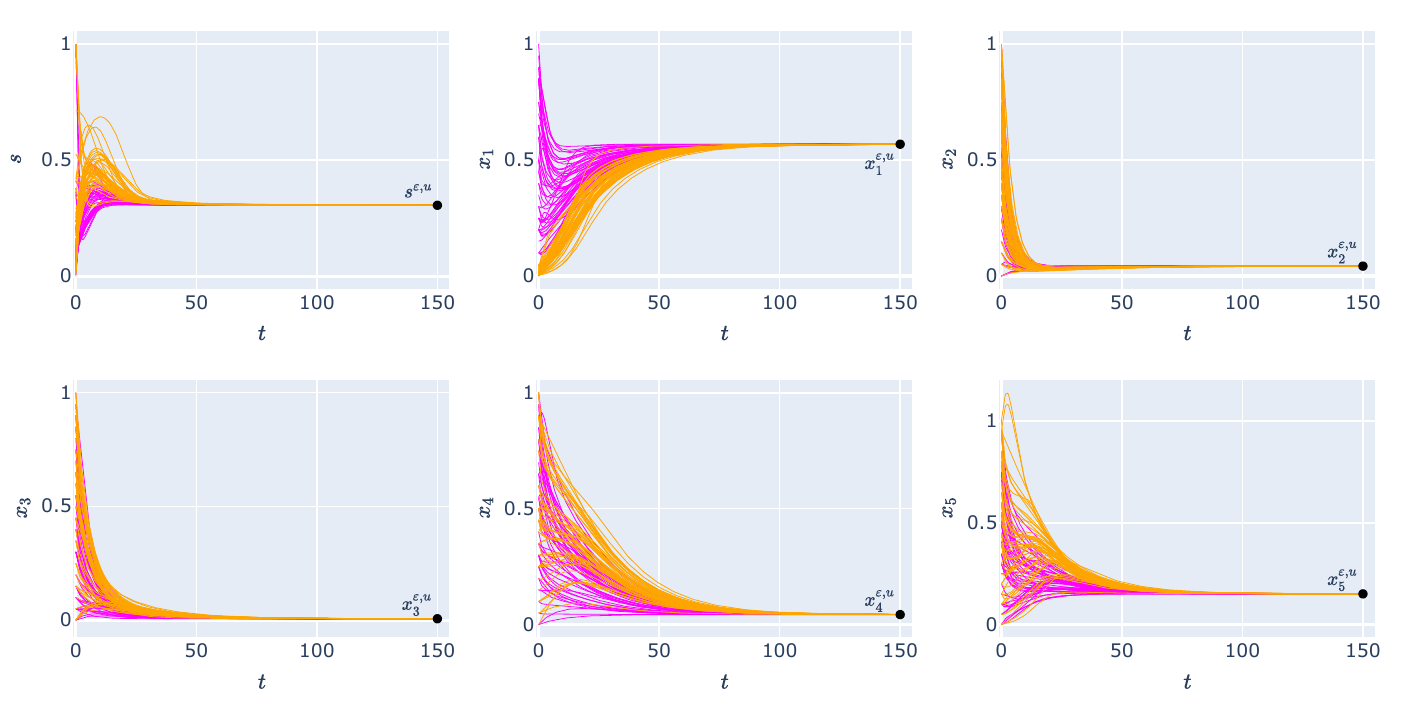}
    \caption{
Solutions to \eqref{eq:cultivated-CS} for 50 randomly generated initial conditions in the set $\Delta_{1,\alpha}$ (magenta) and 50 others in the set 
$\Delta\backslash \Delta_{1,\alpha}$ (orange). }
    \label{fig:monod_cul2}
\end{figure}

\begin{figure}[H]
    \centering
    \begin{subfigure}[b]{0.45\textwidth}
        \centering
        \includegraphics[width=1.1\linewidth]{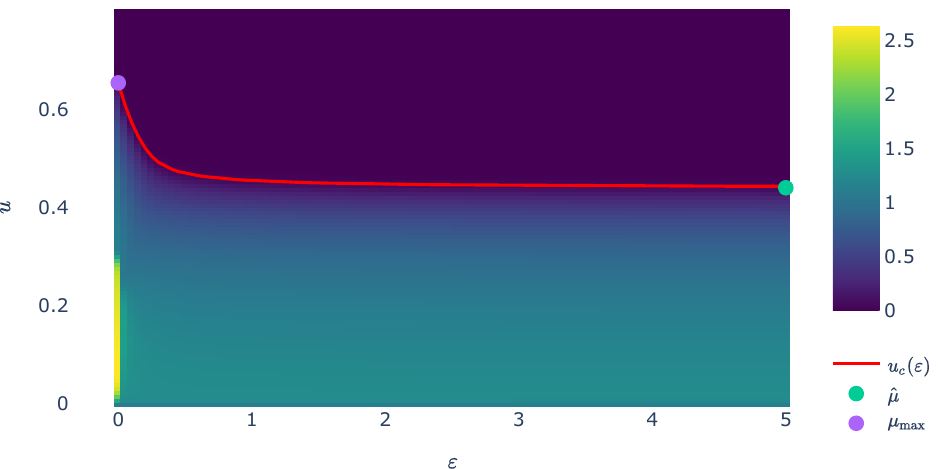}
        \caption{Case 1: constant perturbation matrix.}
        \label{fig:d0_lineal}
    \end{subfigure}
    \hfill 
    \begin{subfigure}[b]{0.45\textwidth}
        \centering
        \includegraphics[width=1.1\linewidth]{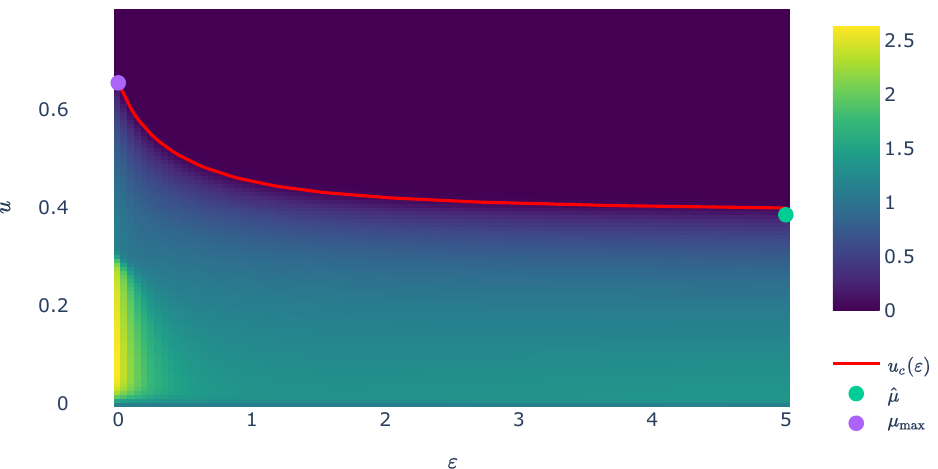}
        \caption{Case 2: non-constant perturbation matrix.}
        \label{fig:d0_cult}
    \end{subfigure}
    \caption{Plot of $\eps \mapsto u_c(\eps)$ in red. The color indicates the distance between the solution of \eqref{eq:lineal-CS} 
    at time $t_f:=200$ and the wash-out for the fixed initial condition {\small{$(0.15,0.15,0.15,0.15,0.15,0.25)$}} and for $(\eps,u)\in [0,10] \times [0,0.8]$.}
    \label{fig:d0}
\end{figure}

\begin{figure}[H]
    \centering
    \begin{subfigure}[b]{0.45\textwidth}
        \centering
        \includegraphics[width=\textwidth]{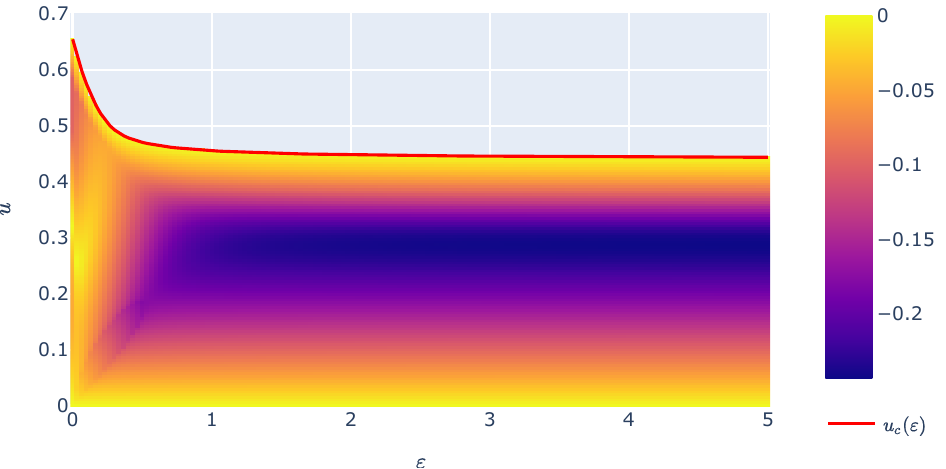}
        \caption{Case 1: constant perturbation matrix.}
        \label{fig:Hurwitz_lineal}
    \end{subfigure}
    \hfill 
    \begin{subfigure}[b]{0.45\textwidth}
        \centering
        \includegraphics[width=\textwidth]{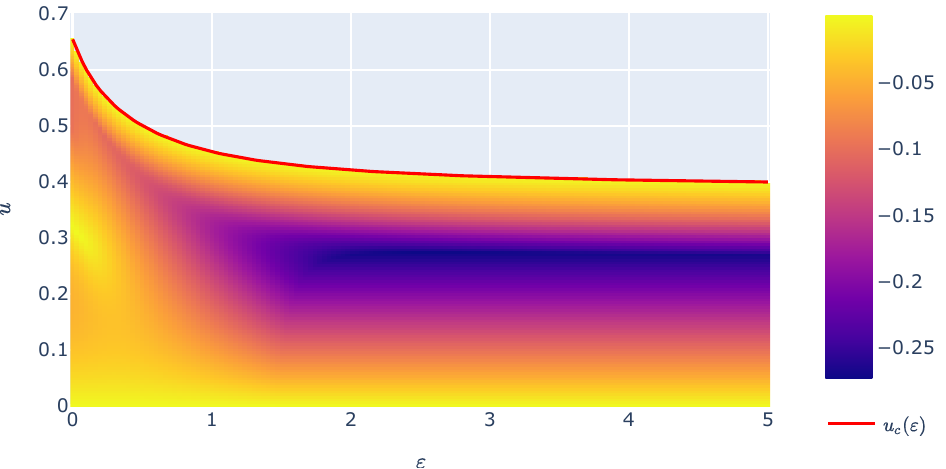}
        \caption{Case 2: non-constant perturbation matrix.}
        \label{fig:Hurwitz_cult}
    \end{subfigure}
    \caption{Plot of $\lambda(J_{\eps,u})$ together with the curve $\eps \mapsto u_{\rm c}(\eps)$ for $(\eps,u)\in [0,5] \times [0,0.7]$.}
    \label{fig:Hurwitz}
\end{figure}

The numerical simulations on systems  \eqref{eq:lineal-CS} and \eqref{eq:cultivated-CS} are intended to highlight their various properties and are organized as follows:
\begin{itemize}
\item[$\bullet$] Fig.~\ref{fig:monod_lin1} and Fig.~\ref{fig:monod_cul1} depict solutions to \eqref{eq:lineal-CS} and to \eqref{eq:cultivated-CS} respectively, for $100$ random initial conditions starting in $(0,10]^n \times [0,s_{in}]$ ; 
\item[$\bullet$] Fig.~\ref{fig:monod_lin2} and Fig.~\ref{fig:monod_cul2} depict solutions to \eqref{eq:lineal-CS} and to \eqref{eq:cultivated-CS} respectively, for $100$ random initial conditions in 
$\Delta_{1,\alpha}$ (magenta) and also in $\Delta \backslash \Delta_{1,\alpha}$ (orange) ; 
\item[$\bullet$] Fig.~\ref{fig:d0_lineal} and Fig.~\ref{fig:d0_cult} represent an operating diagram (in the parameters $\eps$ and $u$) that depicts the separation of the space $(\eps,u)$ 
into two distinct regions via the curve $u_c(\eps)$ ;

\item[$\bullet$] Fig.~\ref{fig:Hurwitz_lineal} and Fig.~\ref{fig:Hurwitz_cult} depict the value of $\lambda(J_{\eps,u})$ in the parameter space $(\eps,u)$, the goal being to check that the Jacobian matrix of \eqref{eq:lineal-CS} at $(x^{\eps,u},s^{\eps,u})$, $J_{\eps,u}$, is  Hurwitz  for 
every $(\eps,u)\in \R_+ \times(0,u_c(\eps))$. 
\end{itemize}

\noindent {\bf{Comments on Fig.~\ref{fig:monod_lin1} \&  Fig.~\ref{fig:monod_cul1}}}. This figure shows that all simulated trajectories converge to the coexistence steady-state $E_{\eps,u}=(x^{\eps,u},s^{\eps,u})$ even if 
initial conditions are taken outside the set $\Delta_{1,\alpha}$. Numerical values of this equilibrium are as follows: 
\begin{itemize}
\item[$\bullet$] for \eqref{eq:lineal-CS}, we find that  $(x^{\eps,u},s^{\eps,u})=(0.48, 0.11,0.04,0.03,0.03,0.37)$ ; 
\item[$\bullet$] for \eqref{eq:cultivated-CS}, we find that $(x^{\eps,u},s^{\eps,u})=(0.57, 0.04,0.006,0.04,0.15,0.31)$. 
\end{itemize}
In contrast with the case $\eps=0$ (corresponding to the chemostat system) where only one species survives, thanks to the CEP, we see that all species are present asymptotically. 
In the first case, dominant species have an index $i$ close to $i^*$, the index of the species that wins the competition for $\eps=0$ (see also~\cite{BCLC}). In the second case, the matrix depends on the substrate, so, this property is not clear as Fig.~\ref{fig:monod_cul1} shows.  

\noindent {\bf{Comments on Fig.~\ref{fig:monod_lin2} \& Fig.~\ref{fig:monod_cul2}}}. As expected, trajectories starting in $\Delta_{1,\alpha}$ converge to the coexistence steady-state, but, this figure  shows that it is also the case if initial conditions are in $\Delta\backslash \Delta_{1,\alpha}$ although it is not predicted by Theorem~\ref{main1}. 
\medskip

\noindent {\bf{Comments on Fig.~\ref{fig:d0}}}. As expected, $\eps \mapsto u_c(\eps)$ decreases from $\mu_{max}$ to $\hat \mu$. 
In both cases, we have $\mu_{max}=0.65$ (see Proposition~\ref{lem-sec3}) whereas in the case of the constant matrix $\Theta$, one has 
$\hat \mu=0.44$ and in the second case, one has $\hat \mu=0.39$. This figure is related to Proposition~\ref{prop:terence} and can be seen as an operating diagram : 
\begin{itemize}
\item[$\bullet$] above the red curve $u_c(\cdot)$, the wash-out is the only steady-state and it is globally asymptotically stable ; 
\item[$\bullet$] on the red curve $u_c(\cdot)$, the wash-out is the only steady-state and it is stable ;
\item[$\bullet$] below the red curve $u_c(\cdot)$, there are two steady-states, namely the wash-out (unstable) and the coexistence steady-state. The behavior of the system toward the coexistence steady-state is described in Theorem~\ref{main1}.  
\end{itemize}
The figure confirms that whenever $u>u_c(\eps)$, then the system converges to the wash-out whereas this is not the case whenever $u<u_c(\eps)$. 
 

\noindent {\bf{Comments on Fig.~\ref{fig:Hurwitz}}}. The figure indicates for every $(\eps,u)\in[0,5] \times [0,0.8]$ such that $u<u_c(\eps)$ the value of
$\lambda(J_{\eps,u})$, the goal being to check numerically that $J_{\eps,u}$ is Hurwitz so that $(x^{\eps,u},s^{\eps,u})$ is locally asymptotically stable. 
Note that $(x^{\eps,u},s^{\eps,u})$ does not exist  for $u>u_c(\eps)$.


\section{Conclusion and perspectives}{\label{sec:conc}}

In this paper, we have addressed the question of global stability of the chemostat system including a perturbation term modeling any interaction type between species. In this context, the perturbation term depends on the species concentration (to model interaction between them) and on a perturbation parameter quantifying the amplitude of the perturbation w.r.t.~the non-perturbed dynamics, but it may also depend on the substrate concentration. We have demonstrated a global stability result toward an equilibrium point in a 
subset $\Delta_{i^*,\alpha}$ of the invariant set (parametrized by a ``small'' parameter $\alpha$) provided that the perturbation parameter $\eps$ is small enough. 
Our methodology relies on the combination of the Malkin-Gorshin Theorem \cite{tewfik0} and a result by Smith and Waltman about the conservation of global attractivity when considering  perturbed dynamical systems \cite{SW99}. The index $i^*$ corresponds to the species that dominates the competition in the chemostat system (without any perturbation). 
Our result may not be sharp, however, as we noted, proving global stability within the invariant set $\Delta$ (let's say, for $\eps$ small enough) is likely a challenging problem. This is because, according to the competitive exclusion principle, solutions to the non-perturbed system do not converge to  the equilibrium point $E_{i^*}$ when the initial condition is on the boundary of $\Delta_{i^*,\alpha}$, specifically when $x_i^*(0)=0$.

Our approach has also shown that one can consider a large class of functions modeling the exchange term although accurate properties of the equilibrium point as given by Theorem~\ref{main1} can be obtained whenever the perturbation term is linear w.r.t.~species (but possibly non-linear with respect to the substrate). It is important to note that the global stability result has extended the study conducted in \cite{BCLC}, which examined stability as the dilution rate approaches zero, to the case where the perturbation parameter tends to zero. 
When the exchange term depends linearly on $x$, a numerical study of the system (asymptotic stability, operational diagram, computation of the Jacobian, etc.) supported the theoretical study conducted in two sub-cases (relying on whether the transition matrix depends on the substrate or not). 
These simulations highlight the fact that as soon as the exchange term is linear, then the globally stable equilibrium is such that all species are present asymptotically. 

It would be interesting to prove that the coexistence steady-state is locally asymptotically stable for all $(\eps,u)\in \R_+^*\times \R_+^*$ such that $u<u_c(\eps)$ (based on a thorough study of the Jacobian matrix which is a rank-one perturbation of a Hurwitz matrix). Additionally, it would be valuable to show the existence of a uniform bound on the perturbation parameter ({\it{i.e.}}, independent of $\alpha$), probably by using another approach (such as Lyapunov functions). This would allow us to obtain a global stability result for all initial conditions in $\Delta$. Another perspective beyond this work would be to carry out a similar analysis but in the case of a 
population dynamics model structured around a single phenotypic trait which can be seen as if we were considering an infinite number of species in the system.

\section*{Acknowledgements}

We are grateful to Tewfik Sari for taking the time to have discussions on global stability of dynamical systems involving regular perturbations and for indicating to us references related to the Malkin-Gorshin Theorem. 

\section{Appendix}

{\it{Proof of Theorem~\ref{prop-claudia-ext-SW}}} (based on \cite[Corollary 2.3]{SW99}). Hereafter, $r(A)$ stands for the spectral radius 
of a matrix $A\in \R^{n \times n}$. 

\smallskip

\noindent {\it{Step 1}}. Let $s>0$ be fixed and recall that $y^*$ is a steady-state of \eqref{eq:generalode} with $\eps=0$, {\it{i.e.}}, $f(y^*,0)=0$. By using standard results about the differentiability of the solution to an ordinary differential equation w.r.t.~the initial condition, we have 
\begin{equation}
    D_{y_0} y(s,y^*,0) = e^{ sD_y f(y^*,0) }\,  \bar 1,
\end{equation}
where $\bar 1:=(1,...,1)\in \R^m$. 
Since $D_y f(y^*,0)$ is a Hurwitz matrix, we deduce that 
$$
r(D_{y_0} y(s,y^*,0)) <1. 
$$
Thus, there exist a norm $|\cdot|$ on $\R^m$ and a real $\rho\in (0,1)$ such that 
$$
\norm{D_{y_0}y(s,y^*,0)}<\rho,
$$
see, {\it{e.g.}}, \cite{zeidler1993nonlinear} for the existence of such a norm. 
By using the continuity of $D_{y_0}y(s,\cdot,\cdot)$, there is $(\bar{\eps}_1,\eta)\in (0,\bar \eps)\times \R_+^*$ (remind that $\bar \eps$ is defined in \ref{hyp:1}) such that\footnote{Here, $B(x,r)$ stands for the open ball in $\R^m$ (for the norm $|\cdot|$) of center $r$ and of radius $r>0$.}  $B(y^*,\eta)\subset U$ and satisfying
        \begin{equation}
        \label{eq:proofsw1}
\forall y_0\in B(y^*,\eta), \ \forall \eps\in [0,\bar{\eps}_1], \;            \norm{D_{y_0}y(s,y_0,\eps)}<\rho. 
        \end{equation}
  By continuity of $y(s,y^*,\cdot)$ w.r.t.~$\eps$, we can also find $\bar{\eps}_2\in(0,\bar{\eps}_1)$ such that
        \begin{equation}
            \label{eq:proofsw2}
 \forall \eps\in[0,\bar{\eps}_2], \;             \norm{y(s,y^*,0)-y(s,y^*,\eps)}<(1-\rho)\eta. 
        \end{equation}
        From \eqref{eq:proofsw1}, we deduce, thanks to the mean value inequality, that 
        \begin{align}
        \label{eq:proofsw3}
  \forall   (y_0,y'_0)\in B(y^*,\eta), \; \forall   \eps\in[0,\bar{\eps}_2], \;      \norm{y(s,y_0,\eps)-y(s,y'_0,\eps)}\leq \rho\norm{y_0-y'_0}. 
        \end{align}
       Combining \eqref{eq:proofsw2} and \eqref{eq:proofsw3} then yields
        \begin{align*}
            \forall   y_0\in B(y^*,\eta), \; \forall   \eps\in[0,\bar{\eps}_2], \;  \norm{y(s,y_0,\eps)-y^*}&\leq \norm{y(s,y_0,\eps)-y(s,y^*,\eps)}+\norm{y(s,y^*,\eps)-y(s,y^*,0)}\\
            &< \rho\norm{y_0-y^*}+(1-\rho)\eta < \eta.
        \end{align*}
        Therefore, for every $\eps\in[0,\bar{\eps}_2]$, $y_0\mapsto y(s,y_0,\eps)$ is a contraction mapping from $B(y^*,\eta)$ into itself. Using the uniform contraction mapping theorem, we deduce that there is a continuous mapping $\hat{y}:[0,\bar{\eps}_2]\to B[y^*,\eta]$ such that $\hat{y}(0)=y^*$ and such that for all $\eps\in [0,\bar{\eps}_2]$, one has $y(s,\hat{y}(\eps),\eps)=\hat{y}(\eps)$. Furthermore, the fixed point $\hat{y}(\eps)$ satisfies the following attractivity property:  
        \begin{equation}
            \label{eq:proofsw4}
\forall (y_0,\eps)\in B(y^*,\eta) \times [0,\bar{\eps}_2], \;         y(ns,y_0,\eps)\overset{n\to \infty}{\longrightarrow} \hat{y}(\eps).  
        \end{equation}
       {\it{Step 2}}. Our objective now is to prove the following: 
               \begin{equation}
            \label{eq:proofsw6}
\exists \bar \eps_3\in (0,\bar\eps_2), \; \forall (y_0,\eps)\in \mathcal{U} \times [0,\bar{\eps}_3], \; 
            y(ns,y_0,\eps)\overset{n\to \infty}{\longrightarrow} \hat{y}(\eps),        
            \end{equation}
      {\it{i.e.}}, that \eqref{eq:proofsw4} holds true for all $y_0\in \mathcal{U}$ (and not only over $B(y^*,\eta)$) up to reducing $\bar \eps_2$. Doing so, we claim the following:  
        \begin{equation}
            \label{eq:proofsw5}
           \exists \bar \eps_3, \;  \forall y_0\in \mathcal{C}, \ \forall \eps\in[0,\bar{\eps}_3], \ \exists m\in \N, \ y(ms,y_0,\eps)\in B(y^*,\eta). 
        \end{equation}
        To prove this property, we proceed by contradiction. It follows that there is a sequence $(\eps_n)_{n\in\N}$ 
        such that $\eps_n\downarrow 0$ and there is a sequence $(y_n)_{n\in\N}$ such that for all $n\in \mathbb{N}$, $y_n\in \mathcal{C}$  and satisfying 
        \[
\forall m\in \mathbb{N}, \; \forall n \in \mathbb{N}^*, \;          \norm{y(ms,y_n,\eps_n)-y^*}\geq \eta.
        \]
        Since $\mathcal{C}$ is compact, we may assume that there is $y_0\in \mathcal{C}$ such that $y_n\to y_0$ (extracting a sub-sequence if necessary). 
        Now, for $\eps=0$, $y^*$ is attracting in $\mathcal{U}$, thus there is $m\in \mathbb{N}$ such that $\norm{y(ms,y_0,0)-y^*}<\frac{\eta}{2}.$
        Moreover, $y(ms,y_n,\eps_n)\to y(ms,y_0,0)$ as $n \rightarrow +\infty$, hence, there is $n \in \mathbb{N}^*$ such that 
        \[
        \norm{y(ms,y_n,\eps_n)-y^*}\leq \norm{y(ms,y_n,\eps_n)-y(ms,y_0,0)} + \norm{y(ms,y_0,0)-y^*} < \eta,
        \]
        which is a contradiction, hence \eqref{eq:proofsw5} holds true. Let now $(y_0,\eps) \in \mathcal{U}\times [0,\bar \eps_3]$. 
        By (H1), there is $T(y_0,\eps)\geq 0$ such that  for all $t\geq T(y_0,\eps)$ one has $y(t,y_0,\eps)\in \mathcal{C}$. Take $m_1\in \mathbb{N}$ such that 
        $m_1s>T(y_0,\eps)$ so that $y(m_1s,y_0,\eps)\in \mathcal{C}$. From  \eqref{eq:proofsw5}, there is $m\in \mathbb{N}$ such that 
        $y((m_1+m)s,y_0,\eps)\in B(y^*,\eta)$. We can now apply \eqref{eq:proofsw4} which implies \eqref{eq:proofsw6} as wanted. 
       
 \smallskip

 \noindent    {\it{Step 3}}. Take $y'_0\in \mathcal{U}$ and let $\eps\in[0,\bar{\eps}_3]$. By (H1) there exists $T(y_0',\eps)\geq 0$ such that 
     for every $t \geq T(y_0',\eps)$, one has $y(t,y'_0,\eps)\in \mathcal{C}$. Thus, we obtain from \eqref{eq:proofsw6}  that 
         \[
         y(t+ns,y'_0,\eps)=y(ns,y(t,y'_0,\eps),\eps)\overset{n\to \infty}{\longrightarrow}  \hat{y}(\eps). 
         \]
 But  \eqref{eq:proofsw6} also implies that $y(ns,y'_0,\eps)\to\hat{y}(\eps)$ as $n \rightarrow +\infty$, hence, since $y(t,\cdot,\eps)$ is continuous, 
         \[
         y(t+ns,y'_0,\eps)=y(t,y(ns,y'_0,\eps),\eps)\overset{n\to \infty}{\longrightarrow} y(t,\hat{y}(\eps),\eps). 
         \]
         We have thus proved that for every $\eps\in[0,\bar{\eps}_3]$ and for every $t\geq T(y_0',\eps)$, one has        
         \begin{equation}
             \label{eq:proofsw7}
           y(t,\hat{y}(\eps),\eps) = \hat{y}(\eps). 
         \end{equation}

We deduce that for $\eps\in[0,\bar{\eps}_3]$, $\hat{y}(\eps)$ is a steady-state of \eqref{eq:generalode}, {\it{i.e.}}, $f(\hat{y}(\eps),\eps)=0$, 
so that by Cauchy-Lipschitz's Theorem, one has 
$y(t,\hat{y}(\eps),\eps) = \hat{y}(\eps)$ for every $t\geq 0$.  
       Finally, let us prove that for every $y_0\in \mathcal{U}$ and every $\eps\in[0,\bar{\eps}_3]$, one has $y(t,y_0,\eps)\rightarrow \hat{y}(\eps)$ as $t\rightarrow +\infty$. By contradiction, suppose that there exist $y_0\in  \mathcal{U}$, $\eps\in[0,\bar{\eps}_3]$, $\nu >0$, and $(t_n)_{n\in\N}$ with  
        $t_n\to +\infty$ such that 
        \[
       \forall n \in \mathbb{N}, \;  \norm{y(t_n,y_0,\eps)-\hat{y}(\eps)} \geq \nu. 
        \]
        Let us write $t_n = k_n s+ \tau_n$ where $k_n\in \N$ and $\tau_n \in [0,s]$. 
        Since $(\tau_n)_{n\in\N}$ is bounded, there is $\tau \in [0,s]$ such that, up to a sub-sequence, $\tau_n\rightarrow \tau\in\big[0,s]$. By using \eqref{eq:proofsw6} and the continuity of $y(\cdot,\cdot,\eps)$, we obtain
        \[
        y(t_n,y_0,\eps)=y(\tau_n, y(k_ns,y_0,\eps),\eps) \overset{n\to\infty}{\longrightarrow} y(\tau, \hat{y}(\eps),\eps)=\hat{y}(\eps),  \]
        where the last equality follows from the fact that $\hat{y}(\eps)$ is a steady-state of \eqref{eq:generalode}. 
         
        Hence, we have a contradiction. 
        We have thus proved that for all $y_0\in \mathcal{U}$ and for all $\eps\in[0,\bar{\eps}_3]$, $y(t,y_0,\eps)\rightarrow \hat{y}(\eps)$ as $t\rightarrow +\infty$ which proves the desired result (taking $\bar \eps_0:=\bar \eps_3$).  \hfill $\Box$
        
\medskip

\noindent {\it{Proof of Theorem \ref{prop-claudia-MG}}}. We recall that the Euclidean norm of $\xi\in \R^m$ is denoted by $|\xi|$. Hereafter, 
$B(x,r)$ stands for the open ball in $\R^m$ (for the Euclidean norm) of center $x$ and of radius $r>0$.
We start by proving the following property:
    \begin{equation}{\label{MK-tmp}}
\forall \delta>0, \; \exists \bar \eps >0, \; \forall \eps\in [0,\bar \eps], \; \forall y_0\in \mathcal{C}, \ \forall t\geq 0, \;       |y(t,y_0,\eps)-y(t,y_0,0)| \leq \delta.
\end{equation}
 
Let then $\delta>0$ be fixed. 
    First, recall that $D_yf(y^*,0)$ is a Hurwitz matrix, hence $y^*$ is exponentially stable for $\eps=0$. Since exponential stability implies B-stability 
    (see  \cite{tewfik0}), Theorem 4 in \cite{tewfik0} implies that there is $\eps_1>0$ such that for all $\eps\in[0,\eps_1]$, one has
        \begin{equation}
        \label{eq:prooftewfik}
 \forall y'\in B(y^*,\eps_1), \; \forall t\geq 0, \;      y(t,y',\eps)\in B(y^*,\tfrac{\delta}{2}).
    \end{equation}
    Observe that $\eps_1\leq\delta/2$, namely, because for any initial condition $y'\in B(y^*,\eps_1)$, one has $y(0,y',\eps_1)=y'\in B(y^*,\delta/2)$ from \eqref{eq:prooftewfik}. 
Second, we claim that $t\mapsto y(t,y_0,0)$ uniformly converges over $\mathcal{C}$ to $y^*$ as $t\rightarrow+\infty$. This result is standard and it follows by combining the local stability property of $y^*$ together with the compactness of $\mathcal{C}$ and the continuity of an ordinary differential equation w.r.t.~initial conditions (that is why, we omit the details). Applying this property then implies 
    \begin{equation}
        \label{eq:prooftewfik2}
\exists T>0, \; \forall t\geq T, \; \forall y_0 \in \mathcal{C}, \;       y(t,y_0,0)\in B(y^*,\tfrac{\eps_1}{2}).
    \end{equation}
   
   Using the continuity of solutions to an ordinary differential equation w.r.t.~a parameter, we get
    \begin{equation}
        \label{eq:prooftewfik3} 
 \exists \eps_2>0, \; \forall \eps\in [0,\eps_2], \; \forall t\in[0,T], \; \forall y_0\in \mathcal{C}, \;      |y(t,y_0,\eps)-y(t,y_0,0)|<\frac{\eps_1}{2}. 
    \end{equation}
   
   Then, combining  \eqref{eq:prooftewfik2} and \eqref{eq:prooftewfik3} yields 
    \[
      \forall \eps\in [0,\eps_2], \;  \forall y_0 \in \mathcal{C}, \;  |y(T,y_0,\eps)-y^*|\leq \eps_1, 
    \]
  hence, for every $y_0 \in \mathcal{C}$ and for every $\eps\in [0,\eps_2]$, one has $y(T,y_0,\eps)\in B(y^*,\eps_1)$. Set $\eps_3:=\min(\eps_1,\eps_2)$ and let us given $\eps\in [0,\eps_3]$ and $y_0 \in \mathcal{C}$.  Applying \eqref{eq:prooftewfik} with $y':=y(T,y_0,\eps)$ yields
  $$
\forall t \geq T, \;   y(t,y_0,\eps)=y(t-T,y(T,y_0,\eps),\eps)\in B(y^*,\tfrac{\delta}{2}).  
  $$
This, together with \eqref{eq:prooftewfik2} implies that 
    \begin{align*}
       \forall t \geq T, \;   |y(t,y_0,\eps)-y(t,y_0,0)| \leq |y(t,y_0,\eps)-y^*|+|y(t,y_0,0)-y^*| 
         \leq \frac{\delta}{2}+\frac{\eps_1}{2} \leq \delta,
    \end{align*}
    since $\eps_1\leq \delta$.  
    Finally, \eqref{eq:prooftewfik3} gives
    $$
    \forall t\in[0,T], \;    |y(t,y_0,\eps)-y(t,y_0,0)|<\delta. 
    $$
   We  can then conclude that
     $$
     \forall \eps\in [0,\eps_3], \; \forall y_0\in \mathcal{C}, \; \forall t \geq 0, \; |y(t,y_0,\eps)-y(t,y_0,0)|<\delta. 
     $$
This ends up the proof of the Theorem setting $\bar \eps:=\eps_3$. 
     
\hfill $\Box$


\end{document}